\documentclass[11pt,twoside]{amsart}%______---IMPOSTAZIONI_RIVISTA
\usepackage{amsmath, amsthm, amscd, amsfonts, amssymb, graphicx, color}
\usepackage[bookmarksnumbered, plainpages]{hyperref}
\newtheorem{thm}{Theorem}
\numberwithin{thm}{section}
\newtheorem{lemma} [thm]{{Lemma}}
\newtheorem{proposition} [thm] {Proposition}
\newtheorem{definition} [thm] {Definition}
\newtheorem{Fact} [thm] {Fact}
\newtheorem{remark}[thm]{Remark}
\newtheorem{example}[thm]{Example}
\newtheorem{cor}[thm]{Corollary}

\newtheorem{Cor}[thm]{Corollary}
\newtheorem{corollary}[thm]{Corollary}
\newtheorem{Prop} [thm] {Proposition}
\newtheorem{problem} [thm] {Problem}
\newtheorem{question} [thm] {Question}
\newtheorem{conjecture} [thm] {Conjecture}

\newcommand{\Fit}{\mbox{\rm Fit}}%FITTING
\newcommand{\End}{\mbox{\rm End}}%ENDOMORFO
\newcommand{\IEnd}{{\mathcal I}\mbox{\rm End}}%ZA
\newcommand{\Hom}{\mbox{\rm Hom}}%ENDOMORFO
\newcommand{\Aut}{\mbox{\rm Aut}}%AUTOMORFO
\newcommand{\PAut}{\mbox{\rm PAut}}%AUTOMORFISMI_POTENZA
\newcommand{\QAut}{\mbox{\rm QAut}}%AUTOMORFISMI_POTENZA
\newcommand{\WAut}{\mbox{\rm WAut}}%AUTOMORFISMI_POTENZA
\newcommand{\IAut}{{\mathcal I}\mbox{\rm Aut}}%ZA
\newcommand{\FAut}{\mbox{\rm FAut}}%AUTOMORFISMI_FINITARI

\newcommand{\St}{\mbox{\rm St}}
\newcommand{\N}{{\Bbb N}}%numeri_naturali_IN_GRASSETTO
\newcommand{\Z}{{\Bbb Z}}%numeri_INTERI
\newcommand{\Q}{{\Bbb Q}}%numeri_rAZIONALI__IN_GRASSETTO
\newcommand{\I}{{\mathcal {I}}}%p-adici_%\cal J
\newcommand{\la}{{\langle}}%parentesi<
\newcommand{\ra}{{\rangle}}%partentesi>
\newcommand{\iso}{\simeq}%accent
\newcommand{\p}{{\varphi}}%VARPHI
\newcommand{\g}{{\gamma}}%
\newcommand{\plus}{\oplus}
\newcommand{\n}{\lhd}%SOTTOGRUPPO_NORMALE
\newcommand{\QED}{\hfill $\square$\medskip}%FINE_DIMOSTRAZIONE_LaTeX2e
\newcommand{\pf}{{\noindent\bf Proof.\ \ }}%PROOF
\newcommand{\im}{\mbox{\rm im}}%IMAGE
\newcommand{\Mod}{\mbox{\rm Mod}}%
\newcommand{\CF}{{\textsc{cf}}}
\newcommand{\BCF}{{\textsc{bcf}}}
\newcommand{\CN}{{\textsc{cn}}}
\newcommand{\BCN}{\textsc{bcn}}

\newcommand{\PP}{\textsc{P}}
\newcommand{\AP}{{\rm AP}}
\newcommand{\BP}{{\rm BP}}
\newcommand{\CP}{{\rm CP}}%
\newcommand{\TP}{{$\tilde{\rm P}$}} 
\newcommand{\sn}{{\rm \kern1pt sn \kern1pt}}%SOTTOGRUPPO_SUBNORMALE
\newcommand{\an}{{\rm \kern2pt an \kern1pt}}%SOTTOGRUPPO_ALMOST_NORMALE
\newcommand{\cf}{{\rm \kern2pt cf \kern1pt}}%SOTTOGRUPPO_con_CORE_FINITO
\newcommand{\f}{{\rm \kern2pt f \kern1pt }}%SOTTOGRUPPO_INDICE_FINITO
\newcommand{\ff}{{\kern1pt \rm \bar{f}\kern1pt }}%SOTTOGRUPPO_INDICE_INVERSO_FINITO
\newcommand{\C}{\rm \kern2pt c \kern1pt}%SOTTOGRUPPO_COMMENSURABILE
\newcommand{\cn}{\rm \kern2pt cn\kern1pt}%SOTTOGRUPPO_COMMENSURABILE_CON_UN_NORMALE
\newcommand{\ssim}{\rm \kern2pt \approx}%STRONGLY_COMMENSURABLE
\def\nbd{neighborhood}
\def\ent{\mbox{\rm ent}\,}
\def\R{\mathbb R}
\def\sF{{\mathcal F}}

\def\nub{\mathrm{nub}}

\begin{document}

%------------------------------------------------------------------------------------%
%%Don not change any thing in this part
%\hskip -0.2 cm
%\begin{tabular}{c r}
%\vspace{-0.6cm}
%%\includegraphics[width=1.9cm]{IJGT-logo}\\
%\href{http://www.theoryofgroups.ir}{\scriptsize  \rm www.theoryofgroups.ir}\\
%\end{tabular}
%\hskip 2 cm
%\begin{tabular}{l}
%\hline
%\vspace{-0.2cm}
%\scriptsize \rm\bf International Journal of Group Theory\\
%\vspace{-0.2cm}
%\scriptsize \rm ISSN (print): 2251-7650, ISSN (on-line): 2251-7669 \\
%\vspace{-0.2cm}
%\scriptsize Vol. {\bf\rm x} No. x {\rm(}201x{\rm)}, pp. xx-xx.\\
%\scriptsize $\copyright$ 201x University of Isfahan\\
%\hline
%\end{tabular}
%\hskip 2 cm
%\begin{tabular}{c c}
%\vspace{-0.1cm}
%%\includegraphics[width=1.9cm]{UI-logo}\\
%\href{http://www.ui.ac.ir}{\scriptsize \rm www.ui.ac.ir}\\
%\vspace{-1cm}
%\end{tabular}
%\vspace{1.3 cm}

%------------------------------------------------------------------------------------%

\title{Inertial properties in groups}
\author{Ulderico Dardano$^*$, Dikran Dikranjan, Silvana Rinauro}

%------------------------------------------------------------------------------------%

\thanks{{\scriptsize
\hskip -0.4 true cm MSC(2010): Primary: 20K30; Secondary: 37A35 (20F28, 16S50, 20E07).
%Received: dd mmmm yyyy, Accepted: dd mmmm yyyy.\\
%$*$Corresponding author
}}
\maketitle
 \vskip -5mm \centerline{\em Dedicated  to the memory of Silvia Lucido}
 \vskip 1mm
%------------------------------------------------------------------------------------%
%This part will be filled in by IJGT
%\begin{center}%Communicated by Patrizia Longobardi
this is a preprint version, the complete paper is going to appear on Int. J. Group Theory
%\end{center}
%------------------------------------------------------------------------------------%

\begin{abstract} 
Let $G$ be a  group and $\p$ be an endomorphism of $G$. A subgroup $H$ of $G$ is called {\em $\p$-inert} if $H^\p\cap H$ has finite index in the image $H^\p$. The subgroups that are {\em $\p$-inert} for all inner automorphisms of $G$ are widely known and studied in the literature, under the name {\em inert} subgroups.

The related notion of {\em inertial endomorphism}, namely an endomorphism $\p$ such that all subgroups of $G$ are  {$\p$-inert},
was introduced in \cite{DR1} and thoroughly studied in  \cite{DR2,DR4}. 
The ``dual" notion of {\em fully inert subgroup}, namely a subgroup that is {\em $\p$-inert} for all endomorphisms of an abelian group $A$,
was introduced in \cite{DGSV} and further studied in  \cite{Ch+, DSZ,GSZ}. 

The goal of this  paper is to  give an overview of up-to-date  known results, as well as some new ones, and show how some applications of the concept of inert subgroup fit in the same picture even if they arise in  different areas of algebra.  We  survey on classical and recent results on groups whose inner automorphisms are inertial. Moreover, we show how 
inert subgroups naturally appear in the realm of locally compact topological groups or locally linearly compact topological vector spaces,
and can be helpful for the computation of the algebraic entropy of continuous endomorphisms.  
\end{abstract}

%\vskip 0.2 true cm

%------------------------------------------------------------------------------------%

\pagestyle{myheadings}
\markboth{\rightline {\sl Int. J. Group Theory x no. x (201x) xx-xx \hskip 5.5 cm
U. Dardano, D. Dikranjan, S. Rinauro }}
         {\leftline{\sl Int. J. Group Theory x no. x (201x) xx-xx \hskip 5.5 cm
U. Dardano, D. Dikranjan, S. Rinauro }}

%%%%%%%%%%%%%%%%%%%%%%%INIZIO_ARTICOLO%%%%%%%%%%%%%%%

% {\small \tableofcontents}

%%%%%%%%%%%%%_SECTION_1_%%%%%%%%%%%%%%%%%%%%%

%\input{input.tex}

\section{Preliminaries}\label{Preliminaries}

\subsection{Introduction}\label{Introduction}

A subgroup $H$ of a group $G$ is  {\em inert} if the index $|H : H^g \cap H|$ is finite for all $g\in G$.  
Clearly normal and finite subgroups are inert. Many  generalizations of the notion of normality have appeared in the literature (see for example \cite{FdG}). 
The concept of inert subgroup allows to put many of 
  them in a common framework. It seems to have been introduced in 1993 in papers of Belyaev \cite{B,B2} as a tool in the investigation of infinite simple groups.  In a survey paper on this subject  \cite{NATO} Belyaev gives credit to Kegel for coining the term {\it inert subgroup}. The concept of inert subgroup tacitly involves inner automorphisms of a group and has been therefore extended to general automorphisms (\cite{DR1}) in 2012. Then in \cite{DGSV} the concept has been  extended to general endomorphisms of abelian groups. Let us give a basic definition which  reveals to be a powerful tool. One may replace the relation of invariance of a subgroup with respect to an endomorphism by a   weaker condition which is trivially satisfied by all finite subgroups.

\begin{definition}\label{def-inertial}   Let  $\p$ be an endomorphism and $H$  a subgroup of a group $G$. Then $H$  is called {\em $\p$-inert} if $H^\p\cap H$ has finite index in the image $H^\p$.
\end{definition} 

A subgroup $H$ is {\it inert} in a group $G$ if and only if $H$ is $\p$-inert in $G$ for each inner automorphism $\p$ of $G$.  
The relevant novelty in the point of view adopted in \cite{DGSV} is imposing the relation from Definition \ref{def-inertial} between a {\em single pair} of a subgroup and an endomorphism,
allowing  one to  extend it further to arbitrary pairs of a family of subgroups and a family of endomorphisms.  In \S\ref{uniformly} we also recall the possibility to distinguish ``uniformly'' inert subgroups as in \cite{BL}.

The concept of inert subgroup has many applications in different areas of algebra. In this paper we try to collect some of them and continue the work done by previous surveys as \cite{NATO,FdG,GSs} even if we are conscious that we cover only a  part of the  subject. Many of these results have been presented at several meetings ``Ischia Group Theory'', an environment that stimulated many interesting and fruitful discussions. Here we also present some new results, propose problems and topics for  further investigations. 

This survey is organized as follows. In the rest of this section (\S\S 1.2--1.4) we give an overview of the paper. In Section \ref{subgroups} we will consider groups in which all or many subgroups are inert or have a stronger property. 
In Section \ref{Groups} we will consider groups of automorphisms $\p$, which have the property that each subgroup is $\p$-inert and $\p^{-1}$-inert (as in \cite{DR1,DR3}) or some stronger property.  
In Section \ref{abelian} we will consider endomorphisms $\p$ of an abelian group $A$ for which each subgroup of $A$ is $\p$-inert and the ring that they form (as in \cite{DR2,DR4}).  
In Section \ref{Fully} we will dually consider  fully inert subgroups, that is subgroups of $A$ which are $\p$-inert with respect to all endomorphisms $\p$ of $A$. 
 There are satisfactory results for many classes of abelian groups  (as in  \cite{DSZ,GSZ}). Moreover, we will describe abelian groups which are inert in their divisible hull (as in \cite{DGSV}). We will also consider some variations of these problems  concerning the relation of fully inert subgroups to fully invariant ones (as in \cite{Cal,GSZ-J}).

 In Section \ref{appl} we recall the algebraic entropy functions $\ent$ and $h_{alg}$ 
for group endomorphisms, as well as the notion of growth of group endomorphism (as in \cite{DG*,DG3,DG:PC,GBSp,GBS+}), which naturally extends the classical function of growth studied by Milnor and Gromov in geometric group theory \cite{dH}.  
Following  \cite{DGSV}, we show how the use of $\p$-inert subgroup leads to a better behaved intrinsic entropy funtion $\widetilde{\ent}$
which works well also in non-torsion abelian groups. Among others, we discuss an extension of the entropy function $\widetilde{\ent}$ to the case of 
non-abelian groups. The final part of the section offers various examples, from the realm of topological groups or topological vector spaces, where 
inert subgroups appear in a natural way and are helpful, again, for the computation of the entropy.  

\subsection*{Notation, terminology and basic facts}
We recall here briefly some notation, definitions and basic facts of those used more frequently in the paper. For the sake of brevity, we omit others that  can be found in the monographs  \cite{F,LR, Rm}.

The set of positive natural numbers will be denoted by $\N$. We use both multiplicative and additive notation for group operation and correspondingly right- and left-hand notation for endomorphisms. As usual, the {\em order} of a group $G$ is  the cardinality of the underling set and the order of an element $g$ of $G$ is the order of the subgroup $\la g\ra$ generated by $g$.
We denote by $id_G$ (or simply $id$, when no confusion is possible) the 
identity map of $G$. If $\p$ is an endomorphism of $G$, we say that a subgroup $H$ of $G$ is $\p$-invariant if $\p(H)\subseteq H$. If $\Gamma$ is a group acting on $G$, then we say that a subgroup $H$ of $G$ is $\Gamma$-invariant if $H^\g\le H$ for each $\g\in\Gamma$. We denote  $H_\Gamma:=\cap_{\g\in\Gamma}H^\g$ (resp. $H^\Gamma:=\la H^\g | {\g\in\Gamma}\ra$) the largest (resp. smallest) $\Gamma$-invariant subgroup of $G$ contained in (resp. containing) $H$.  In particular, $H_G$ is the core of $H$ in $G$ and $H^G$ is the normal closure of $H$ in $G$.

We use $A$ to denote exclusively abelian groups, while $t(A)$, $A_\pi$ and $D(A)$ denote the torsion subgroup,
the $\pi$-component  and the maximum divisible subgroup of $A$, respectively.  Further,  if $n\in\N$, we denote  $A[n]:=\{a\in A\ |\ na=0\ \}$. If $nA=0$, we say that $A$ is \emph{bounded} by $n$ and the least such $n$ is called the exponent of A. We just say that $A$ is bounded or that $A$ has finite exponent, if $A$ is bounded by some $n$. A  Pr{\"u}fer (or {\em quasi-cyclic}) $p$-group is a group isomorphic to the $p$-component $\Z(p^\infty)$ of the quotient $\Q/\Z$ of the additive group of the rational numbers.  If $n\in\N$ then $\pi(n)$, is the set of prime divisors of $n$. We denote $\Q^n:=\Q\oplus\ldots\oplus\Q$ ($n$ times).  If $\pi$ is a  set of primes, we denote by $\pi'$ the set of primes not in $\pi$ and by $\Q^{\pi}$ the ring (or just the additive group) of rational numbers whose denominator is divisible by primes in $\pi$ only. 
 
 Recall that  the torsion free rank $r_0(A)$ of an abelian group $A$ is the cardinality of any maximal  independent subset of elements of infinite order.  Then $r_0(A)<\infty $ if and only if $A/t(A)$ is isomorphic to a subgroup of $\Q^n$ for some $n\in\N$. If $p$ is a prime, $r_p(A)$ denotes  the cardinality of any maximal independent subset of elements of $A$ with order $p$. It follows that $r_p(A)<\infty$ if and only if $A_p$ is the direct sum of finitely many cyclic or quasi-cyclic groups. Furthermore, the rank of $A$ is $r(A):= \max\{r_0(A), \sup \{r_p(A): p \ prime\}\}$.
 
If $\mathfrak{P}$ and  ${\mathfrak{P}}_1$ are properties for groups, then we say that a group $G$ is {\em locally}  $\mathfrak{P}$ if each finitely generated subgroup of $G$ has  $\mathfrak{P}$ and that $G$ is   $\mathfrak{P}$-{\em by}-$\mathfrak{P}_1$ if $G$ has a normal subgroup $N$ such that $N$ has  $\mathfrak{P}$ and $G/N$ has  $\mathfrak{P}_1$. A group $G$ is said  {\em metabelian} if its derived subgroup $G'$  is abelian. Moreover, 
a group is called {\em locally  graded} if each of its non-trivial finitely generated subgroups  has a non-trivial finite image.  
 A group is {\em residually} $\mathfrak{P}$ if the trivial subgroup is the intersection of  the set of normal subgroups $N$ such that $G/N$ has $\mathfrak{P}$. Recall that the class of locally graded groups contains all locally (soluble-by-finite) groups and all locally residually finite groups. 

Recall that a subgroup $H$ of a group $G$ is said {\em subnormal} in $G$ if there is a finite series $H=H_0\n H_1\n\ldots\n H_n= G$ where each term is normal in the successive one. A {\em Chernikov group}  is a group which has a  subgroup of finite index which is a direct product of finitely many  Pr{\"u}fer groups (for possibly different primes).

\vskip-1cm
%%%%%%%%%%%%%%%%%%%%%%%%%%%%%%%%%%%

\subsection{Inert subgroups}\label{Inert} 
Let us introduce some further terminology.

\vskip-1cm
\begin{definition}  If $H$, $K$ are subgroups of a group $G$ we say that:
\begin{enumerate}
\item $H$ is {\em almost contained} in $K$ (and write $H\le_a K$) if  $|H:(K\cap H)|<\infty$;
\item $H$ is {\em commensurable with} $K$ (and write $H\sim K$) if $H\le_a K$ and $K\le_aH$.
\end{enumerate}
\end{definition} 

In these terms, for an endomorphism  $\p$ of $G$, a subgroup $H$ is $\p$-inert  if and only if the image $H^\p$ is almost contained in $H$. The relation $\le_a$  is reflexive and transitive;  the equivalence relation that it generates is usually called 
{ commensurability}. Two subgroups $H$ and $K$ are {\em commensurable} if their intersection $H\cap K$ has finite index in both of them, that is $H\le_a K$ and $K\le_a H$. Thus, a subgroup of a group is  inert if and only if it is {commensurable } with each of its conjugates.  

Accordingly, in a recent paper  \cite{PW} in the setting of Model Theory, Placin and Wagner consider the {\em almost normalizer} $\tilde N_G(H)$ as the subgroup formed by the elements $g\in G$ such that $H^g$ is commensurable with $H$.

Commensurability is compatible with taking the intersection of finitely many subgroups and taking quotients. It is preserved under projectivities, i.e. subgroup lattice isomorphisms, even when not induced by group isomorphisms, as proved by Zacher \cite{Z} in 1980. Commensurability appears 
 explicitly mentioned  in 1987  by Heineken and Specht  \cite{HS} who studied groups with only finitely many commensurability classes. Later, in 1997, Heineken and Kurdachenko in \cite{HK} considered groups whose {\em normal} subgroups are partitioned in finitely many commensurability classes. 

Not surprisingly, the structure of a proper  inert subgroup in a simple group is far from arbitrary. For example,  it must be residually finite. We do not go deeper in this direction leading somewhat far from the scope of this survey (see  \cite{B,B2, NATO,DES5,DES6} for more details).

Clearly, a subgroup commensurable with an inert  subgroup, is inert itself. In particular, subgroups which are commensurable with a normal subgroup are inert. Therefore, a subgroup $H$ of a group $G$ having a normal core $H_G$ of finite index $|H:H_G|$, is inert. Such subgroups are called {\em normal-by-finite} in  \cite{BLNSW,CKLRS, FdG}.  
Dually, a subgroup $H$ of a group $G$ having finite index $|H^G:H|$ in its normal closure $H^G$ is inert.  Such subgroups are called {\em nearly-normal} in  \cite{FdG}.

Nevertheless, an inert subgroup need not be  commensurable with a normal subgroup, as shown by the following example due to Peter Neumann. Let $G$ be the  alternating group of an arbitrary infinite set $S$, i.e. the group of  permutations of $S$ which move only finitely many elements and are product of an even number of transposition. Let $P$ be any partition of $S$ into finite subsets, each of cardinality $>2$, and $H$ the subgroup of elements of $G$ which carry each member of $P$ into itself and act on each of these as an even permutation. Then $H$ is inert but not commensurable with a normal subgroup of $G$ (as $G$ is simple), while $H$ has both order and index equal to $|S|$. Other non-trivial examples of inert subgroups are permutable subgroups (see \cite{R}). Further examples and basic properties of inert subgroups can be found in  \cite{B,R}.  Notice that all the subgroups of the so-called  dihedral group (over $\Z$)  $D:=\la x,y \ |\  x^2=1, y^x=y^{-1}\ra=\la y\ra\rtimes \la x\ra$ are inert.

%%%%%%%%%%%%%%%%%%%%%%%%%

\subsection{The inertial  correspondence}\label{basic} 
The straightforward proof of the following fact is omitted.
 
\begin{Fact}\label{Fact} In the lattice of subgroups of a group, the relation  $\sim$ of commensurability is the equivalence generated by the transitive  relation $\le_a$. Both relations are compatible with the meet operation (i.e. intersection) and endomorphisms. 

If the underlying group $A$ is abelian, then the relations $\le_a$ and  $\sim$ are compatible with  
the join operation (i.e. sum) as well.  Thus,  
commesurability is a  congruence in the lattice $l(A)$ of subgroups  of $A$.\qed 
\end{Fact}

Consequently, {the intersection of finitely many inert subgroups is inert. However the intersection of infinitely many inert subgroups need not be inert} (see Corollary \ref{intersezione_inerti}). Further, the join of a pair of inert subgroups of a group $G$ need not be inert, if $G$ is non-abelian. This can be seen as in \cite{R} by considering the soluble group 
$$
G:=\la a,b,x,y \ |\ a^b=a, x^2=1, y^x=y^{-1}, a^x=a^{-1}, b^x=b,  a^y=ab, b^y=b\ra=\la a, b\ra\rtimes \la x,y\ra,
$$
 where $A:=\la a, b\ra$ is free-abelian of rank $2$ and $D:=\la x,y\ra$ is a dihedral group. The infinite subgroup $D$ is generated by the finite subgroups 
$\la x\ra$ and $\la xy\ra$, but $D$ is not inert in $G$, as $D\cap D^a=1$. 

Let us now introduce the dynamical aspect of the inertial relation given by Definition \ref{def-inertial}. Note that, when the underlying group is abelian (in additive notation), a subgroup $H$ is $\p$-inert  if and only if the factor group $(H+\p(H))/H$  is finite.

\begin{Fact}\label{Fact+} Let $A$ be an abelian group, $\Phi$ a subset of the ring $End(A)$ of endomorphisms of $A$ and  $\mathcal{H}$ a subset of the lattice  $l(A)$ of subgroups of $A$. Then: 
\begin{enumerate}
\item the set $\mathcal I_\Phi(A)$ of subgroups which are  $\p$-inert for each $\p\in\Phi$ is a sublattice of $l(A)$ closed under commensurability;
\item the set $\mathcal I_\mathcal{H}\End(A)$  of endomorphisms $\p$ for which each $H\in \mathcal{H}$ is $\p$-inert  is a subring of $End(A)$.\qed
\end{enumerate}
\end{Fact}

\pf The proof is rather straightforward and follows from \cite{DR2,DGSV}. Here we only recall an argument that shows that $\mathcal I_\Phi(A)$ is closed under commensurability. If $A$ is additively written and $H_1\sim H \in \mathcal I_\Phi(A)$, then for every endomorphism $\p\in \Phi$  
each of the factors in the following series is finite 
$$ 
H\cap H_1 \le H  \le H+\p(H) \le H+H_1+\p(H) \le H+H_1+\p(H)+\p(H_1).
$$
Hence, the same holds for the intermediate factor in the series  
$$
H\cap H_1 \le H_1\le H_1+\p(H_1)\le H+H_1+\p(H)+\p(H_1).
$$
 Therefore, $H_1\in \mathcal I_\Phi(A)$ .
\qed\
\medskip 

 In the sequel we simply write $\mathcal I_{\varphi}(A)$, when $\Phi = \{\varphi\}$ is a singleton, and $\mathcal I\End(A)$ for $\mathcal I_{l(A)}\End(A)$.
This case, introduced in \cite{DR2}, deserves a special attention and will be largely discussed in following sections.   
Note that in the next definition the underlying group $G$ is not necessarily abelian.

\begin{definition}\label{def:inertial:endo}\cite{DR2}
An endomorphism $\p$  of a group $G$ such that $H$ is $\p$-inert for each  subgroup $H$ of $G$ is called an {\em inertial endomorphism}. 
\end{definition}

In other words, an endomorphism $\p$ of a group $G$ is inertial if the image of $H$ is almost contained in $H$ for each subgroup $H$ of $G$. According to Fact \ref{Fact+}(2), the inertial endomorphisms of any abelian group $A$ form a subring, say  $\IEnd(A)$, of the whole ring $\End(A)$. 

 The subring $\IEnd(A)$ contains the ideal $F(A)$ formed by the endomorphisms of $A$ with finite image. We come back to this ring in \S\ref{Ring},  while in  \S \ref{IAut(A)} we focus our attention on the group $\IAut(A)$ generated by the  automorphisms $\g$  for which each subgroup  is $\g$-inert. This group is somewhat larger than the group of invertibles of $\IEnd(A)$, that we will call {\it bi-inertial automorphisms}. Note that the automorphism $\p:x\mapsto 2x$ of $\Q^\N$ is inertial but not bi-inertial, as $\Z^\N$ is not $\p^{-1}$-inert).

The counterpart of Definition \ref{def:inertial:endo} is the following concept proposed in \cite{DGSV}: 
 
\begin{definition} \cite{DGSV}
A subgroup $H$ of an abelian group $G$ is called {\em  fully inert} if $H$ is $\p$-inert for each endomorphism $\p$ of $G$.
\end{definition}

\noindent Clearly, a subgroup which is commensurable with a fully invariant subgroup is fully inert. On the other hand, $\Z$ is a fully inert subgroup of $\Q$, but it is not commensurable to any fully invariant subgroup of $\Q$. In Section \ref{Fully} we will investigate the difference between these two properties and 
will discuss results on fully inert subgroups of (divisible) abelian groups.

%%%%%%%%%%%%%%%%%%%%%%%%%%%%%%%%%%%%%

\subsection{Uniformly inert subgroups}\label{uniformly}

The following definition is motivated by a nice result by  Bergman and Lenstra  \cite{BL} in 1989 (see Theorem \ref{BL}). 

\begin{definition}\label{boundedly} A subgroup $H$ of a group $G$ is called {\em uniformly inert} with respect to a subgroup $K$ of $G$ if there is $n\in\N$ such that $|H: H^g\cap H|\le n $ for each $g\in K$.\\
 If it is necessary to emphasize $n$, we may also say {\em $n$-uniformly inert}. \\ 
 In case $K = G$, we simply say $H$ is {\em uniformly inert} (resp., {\em $n$-uniformly inert}). 
\end{definition}

Note that, since $K$ is a subgroup, in the above definition  we could equivalently write that both indeces $|H: H^g\cap H| $ and $|H^g: H^g\cap H|$ are bounded by some natural number for each $g\in K$, since $|H: H^\g\cap H|=|H^{\gˆ{-1}}: H^{\gˆ{-1}}\cap H|$ for any automorphism $\g$ of $G$.

Roughly speaking, the next theorem states that a subgroup of a group is commensurable with a normal subgroup if and only if it is uniformly  inert. 

\begin{thm}\label{BL}\cite{BL} If $H$ and $K$ are subgroups of a group $G$ then the following
conditions are equivalent:
\begin{enumerate}
 \item $H$ is uniformly inert with respect to $K$; 
 \item $H$ is commensurable with a subgroup  $N$ of $G$ which is normalized by $K$. 
\end{enumerate}
Moreover, in case $H$ is $n$-uniformly inert w.r.t. $K$, then there is $N\le G$ such that $N$ is normalized by $K$, $|N:H\cap N|\le n$, $|H:H\cap N|\le n^{n^n}$ 
and there is also  $M\le G$ such that $N\le M$, $M$ is normalized by $K$,  $M/N$ is finite and $|H:H\cap M|\le n$.
\end{thm}

Recall that the bound $ n^{n^n}$ is due to P. M. Neumann and the particular case when $G=K$ had been previously considered by Schlichting \cite{Schlichting} in 1980 when generalizing a result by Baer \cite{Baer}.  
We state a corollary in additive notation.

\begin{cor}\label{BLcorollario}
 For a subgroup $H$ of  an abelian group $A$  the following are equivalent:
\begin{enumerate} 
\item there is $n$ such that $|H+ \g(H)/H|\le n $ for each $\g\in Aut(A)$;
\item $H$ is commensurable with a characteristic subgroup $K$ of $A$.
\end{enumerate}
Moreover, there is  a characteristic subgroup $K$ of $A$ such that $|(H+K)/(H\cap K)|\le n^{n^n+1}$.
\end{cor}

\pf   It is clear that ($2$) implies ($1$). Assume ($1$) and let $G=A\rtimes \Aut(A)$ be the holomorph group of $A$ and note  that for each $g\in G$ one has
 $$
 |H/(H^g\cap H)|=|H^{g^{-1}}/(H\cap  H^{g^{-1}})|= |\langle H,  H^{g^{-1}}\rangle /H|\le n.
 $$
Thus one may apply Theorem \ref{BL} with $K:=G$. So $H$ is commensurable with a normal subgroup $N$ of $G$. Then $N\cap A$ is a characteristic subgroup of $A$ commensurable with $H$ .
\QED
 
Theorem \ref{BL} and Corollary \ref{BLcorollario} suggest the following definitions. For every abelian group $A$ we denote by $\mathcal Inv(A)$ the sublattice of $l(A)$ consisting of all fully invariant subgroups of $A$, and by 
$\mathcal {I} (A)=\mathcal {I}_{End(A)} (A)$ 
the family of all  fully inert subgroups of $A$. According to Fact \ref{Fact+}, $\mathcal {I} (A)$ is a sublattice of $l(A)$, stable under $\sim$, i.e., a subgroup which is commensurable with some fully inert subgroup of $A$ is also fully inert. Easy examples show that in general $\mathcal Inv(A)$ is not stable under $\sim$, in other words, the family $\mathcal Inv^\sim(A)$ of all subgroups of $A$ commensurable with some fully inert subgroup of $A$ properly contains the lattice $\mathcal Inv(A)$, and it is a lattice on its own, according to Fact \ref{Fact} and the fact that $\mathcal Inv(A)$ itself is a (complete) lattice. 

\begin{definition} 
A subgroup $H$ of an abelian group $A$ is {\em uniformly fully inert}, if there is $n$ such that $| (H+\p(H))/H|\le n $ for each $\p\in End(A)$.
\end{definition} 

By an argument similar to that used in Fact \ref{Fact+} one proves the following: 

\begin{proposition}\label{u-lattice} If $A$ is an abelian group, then 
the set $ \mathcal I_u (A)$ of uniformly fully inert subgroups of $A$ 
 is a sublattice of\ \  $l(A)$ which is closed under $\sim$ and
$$
\mathcal Inv(A)\ \subseteq \ \mathcal Inv^\sim(A)\ \subseteq \ \mathcal I_u(A)\ \subseteq \ \mathcal I(A). \eqno(*)
$$ 
\end{proposition}

 Section \ref{Fully} is mainly dedicated to the following: 

 \begin{problem}\label{3reticoli} Study the interrelations between the above lattices. 
\end{problem} 

%%%%%%%%%%%%%%%%%%%%%%%%%%%%%%%%%%%%%%%%%%%

%%%%%%%%%%%%%%____SECTION_2____%%%%%%%%%%%%%%%%%%%%%%%%
%%%%%%%%%%%%%%%%%%%%%%%%%%%%%%%%%%%%
%\medskip 
\section{Groups with many inert subgroups} \label{subgroups}

 It is well known that Dedekind groups, that is groups $G$ in which all subgroups are normal, have a restricted structure (see \cite{R}). In particular, if $G$ is a non-abelian Dedekind group, then the groups $G/Z(G)$ and $G'$ have respectively order $4$ and $2$.

%%%%%%%%%%%%%%

\subsection{Groups in which each subgroup is inert}\label{TIN}
A group is called  {\em inertial}  (or totally inert, or TIN) if all its subgroups are inert. Clearly, all finite groups and Dedekind groups are inertial.
Moreover, if a group $G$ has a finite normal subgroup $F$ such that $G/F$ is an  inertial group, then $G$ is inertial. 

The first result we recall on inertial groups is a theorem by Beidelman and Heineken \cite{BH} in 2000 on groups which are inertial in a rather trivial way.

\begin{thm}\cite{BH} Let $G$ be an infinite group. Then  all infinite 
subgroups  of $G$ are normal if and only if one of the
following holds:
\begin{enumerate}
\item $G$ is a Dedekind group; 
\item  $G$ is an extension of a Pr{\"u}fer $p$-group by a finite Dedekind group;  
\item   $G$ is an extension of a non-abelian infinite group $M$ by a Dedekind group and all infinite  subgroups of $G$ contain $M$.
\end{enumerate}
\end{thm}

 We state now a result by Robinson  \cite{R} in 2006 which guarantees that the class of inertial groups is rather large.
Recall that a group $G$ is said to be an {\em FC-group} if each element of $G$ has a finite number of conjugates. Note that 
the number of conjugates of an element $a$ of $G$ equals the index of the centralizer $C_G(a):=\{g\in G\ | \ ag=ga\}$.

\begin{thm} \cite{R} FC-groups are inertial. 
\end{thm}

The class of inertial groups can be very complex. For example all {\em quasi-finite} groups (that is groups in which all proper subgroups are finite) are inertial. In particular 
 Tarski $p$-groups (that is,  infinite simple groups in which all proper non-trivial subgroup have  order $p$) are inertial. 
 On the other hand,  Belyaev, Kuzucuoglu and  Seckin \cite{BKS} in 1999 have proved:

\begin{thm} \cite{BKS}\label{BKS} { \it No infinite locally finite simple
group can be inertial}.
\end{thm}

Dixon,  Evans and Tortora \cite{DET} in 2010 have extended this result to simple locally graded groups. 

\begin{thm}\cite{DET}
Let $G$ be an inertial infinite simple group. Then $G $ is a finitely generated periodic group without involutions such that $C_G(x)$ is finite for all $1\ne x \in G$ and in which every proper subgroup is residually finite. In particular, $G$ is not locally graded.
\end{thm}

\begin{cor}\cite{DET}
An inertial simple group $G$ of finite exponent 
is quasi-finite. 
\end{cor}

It seems to be unknown whether any inertial simple group is necessarily quasi-finite. 
By restricting the investigation to the realm of generalized soluble groups, Robinson  \cite{R}  in 2006 has proved the following theorem. Recall that a
hyper-(abelian-by-finite) group  is a group which has an ascending (i.e. well ordered) series of normal subgroups whose factors are either abelian or finite. We already noticed that a dihedral group over $\Z$ is inertial.

\begin{thm} \cite{R}\label{RobB}  Let $G$ be a finitely generated hyper-(abelian-by-finite) group. Then $G$ is
 inertial if and only if $G$ has a torsion-free abelian normal subgroup $A$ of finite index in which
elements of $G$ induce by conjugation either the identity or the inversion map.
\end{thm}

From \cite{R} we recall that in a polycyclic group every subgroup is the intersection of subgroups of finite index, but not every subgroup of a polycyclic group is necessarily inert.   According to Theorem \ref{RobB}, this is witnessed   by the free nilpotent group of class and rank $2$. Thus, we have:

\begin{Cor} \label{intersezione_inerti}The intersection of an infinite family of inert subgroups need not to be inert.
\end{Cor}

It is easy to show that if a group $G$ has normal subgroups $F$ and $A$, with $F \le A$, such that $|F|$ and $|G / A|$ are finite, $A/F$ is abelian and each subgroup of $A/F$ is normalized by $G$, then the group $G$ is inertial.   In such a case, one says that $G$ is {\em inertial of elementary type}.  

Concerning non finitely-generated inertial groups, from \cite{LR} we recall that a soluble group $G$ has {\em finite abelian total rank}, or is a soluble FATR-group, if $G$ has a finite series in which each factor is abelian with finite rank and is either torsion-free or a $p$-group.  In 
\cite{R}  D. Robinson has described soluble inertial FATR-groups. We remark that also groups of  non-elementary type  occur, but here we do not report  the classification as it requires  technical details.  We just deduce the following.

\begin{thm} \label{FATR}\cite{R} Let $G$ be a soluble inertial FATR-group. Then $G$ is finite-by-metabelian-by-finite.
\end{thm}

\subsection{Groups in which each subgroup is strongly inert}\label{strongly-inertial} 
Another way to avoid dealing with quasi-finite groups 
is strengthening the condition of being inert. Thus, De Falco, de Giovanni,  Musella,  Trabelsi  \cite{DGMT} in 2013 have considered the following definition.

\begin{definition}\cite{DGMT} A subgroup $H$ of a group $G$ is said to be {\em strongly inert}
if $H$ has finite index in $\langle H,H^g\rangle$ for each $g\in G$.

A group in which all subgroups are strongly inert is called {\em strongly inertial}.
\end{definition}

 Clearly, every { nearly-normal} subgroup is strongly inert. Moreover, if the commutator subgroup $G'$ of
$G$ is finite of size $n$, then the index of $H$ in $\langle H,H^g\rangle$ is at most $n$ for each $g\in G$ and $H$ is strongly inert. 
 Then we recall a celebrated theorem by Neumann \cite{BHN} in 1955 which deals in fact with strongly inertial groups.

\begin{thm}  \cite{BHN}   In a group $G$ each subgroup  is nearly-normal 
 if and only if the group has finite derived subgroup (i.e. $G$ is finite-by-abelian).
\end{thm} 

 Clearly permutable subgroups are strongly inert (see \cite{DGMT}).  Other examples of  strongly inert subgroups are inert  subgroups which are subnormal (see \cite{DR5}). On the other hand, finite subgroups in a dihedral group (over $\Z$) are not strongly inert (although trivially inert). We also recall that in an FC-group all subgroups are  strongly inert (see \cite{DGMT}).

A finitely generated group is strongly inertial if and only it it is an  FC-group, as we see in next theorem.

\begin{thm} \cite{DGMT}  Every finitely generated strongly inertial group is centre-by-finite.
\end{thm}

\begin{cor} Let $G$ be a group whose finitely generated subgroups are strongly inert. Then the commutator subgroup $G'$ of $G$ is locally finite, i.e. $G$ is (locally finite)-by-abelian.
\end{cor}

Thus, a strongly inertial simple group must be finite, since otherwise  $G'=G$ is locally finite, contradicting Theorem \ref{BKS}.

 Recall that a group $G $ is called {\em minimax} if it has a series of finite length whose factors satisfy either the minimal or the maximal condition on subgroups. Recall also that  if $G$ is any soluble-by-finite minimax group, then its finite residual $R$ (i.e. the intersection of all subgroups of finite index) is the direct product of finitely many  Pr{\"u}fer subgroups.  If this subgroup $R$ is actually a $p$-group we say that $G$ is {\em $p$-primary}. Then in  \cite{DGMT} a full classification of soluble-by-finite $p$-primary minimax strongly inertial groups is obtained.
We omit  details, but we recall a relevant consequence.

\begin{cor} \cite{DGMT} A soluble-by-finite minimax group is strongly inertial if and only if it is inertial and has no infinite dihedral sections.
\end{cor}

%%%%%%%%%%%%%%%%%%%%%%%%%%%%%%%%

\subsection{Groups in which each subgroup is normal-by-finite}\label{normal-by-finite} 
We are going to introduce a condition dual to nearly-normality. Recall that a subgroup of a group $G$ is called {\em normal-by-finite} (or sometimes {\em almost-normal}) if 
$H/H_G$ is finite. A normal-by-finite subgroup is inert.

\begin{definition}\label{def-CF} Let $G$ be a group.\hfill  
\begin{enumerate}
\item $G$  is a {\em $CF$-group} if  $\forall H\le G\ \ |H/H_G|< \infty$;
\item $G$ is a {\em BCF-group} if $\exists n\in \N\ \  \forall H\le G\ \ |H/H_G|\le n$.
\end{enumerate}
\end{definition} 
 Since normal-by-finite subgroups are also called ``core-finite" by some authors,  
here CF stands for ``core-finite" and BCF for ``boundedly core-finite". The class of CF-groups seems to be difficult to handle since all quasi-finite groups are CF. In fact normal-by-finite (and even finite) subgroups are not necessarily strongly inert. However, the problem is tractable in a restricted class of groups. It has been treated by   Buckley, Cutolo, Lennox, Neumann,  Rinauro,  Smith and Wiegold around 1995 in \cite{BLNSW, CKLRS, SW}. 

\begin{thm} \cite{BLNSW,SW}
A CF-group $G$ whose periodic quotients are locally finite is abelian-by-finite.  Moreover, if $G$ is periodic, then $G$ is BCF. 
\end{thm} 

Note that every group which has an ascending normal series whose factors are either locally nilpotent or finite has  the property that  periodic quotients are locally finite.

\begin{problem} Is every locally graded  CF-group abelian-by-finite?
\end{problem} 

 Let us see how  periodic  abelian-by-finite CF-groups may be described in a satisfactory way using also results from
  \cite{DR1, FGN}. Recall that an automorphism $\gamma$ of a group $G$ is said a {\it power automorphism} if $H^\g=H$ for each subgroup $H\le G$. \marginpar{} It is well-known  (see \cite{Rm}) that, if $G$ is an 
 abelian $p$-group, then there exists a $p$-adic integer $\alpha$ such that $a^\g=a^\alpha$ for all $a\in G$. Here $a^\alpha$ stands for $a^n$, where $n$ is any integer congruent to $\alpha$ modulo the order of $a$. 

\begin{thm}
A {locally finite group $G$ is CF if and only if} it has a normal abelian subgroup \ $A=C\times (D\times E)$\  of finite index and such that: 
\begin{enumerate}
\item$\pi(D)=\pi(E)$ is finite and $\pi(C)\cap\pi(DE)=\emptyset$; 
\item$D$ is divisible of finite rank $r$ and $E$ has finite exponent $e$; 
\item$G$ induces by conjugation power automorphisms  on $C$, $D$, $E$.
\end{enumerate}
For such a group $G$ one has $|H/H_G|\le e^r\cdot |G/A|$ for every $H\le G$. Therefore, $G$ is a BCF-group.
\end{thm}

Concerning the non-periodic case, recall that a power automorphism of a non-periodic abelian group is either the identity or the inversion map, say for short that it is {\em the power $\pm1$}  respectively. Consider the group  $G:= A \rtimes   \la\g\ra$ where $A=C_\infty\times C_{p^\infty}$
(here $C_\infty$ and $C_{p^\infty}$ are an infinite cyclic group and a Pr{\"u}fer $p$-group)  
and $\g$ is the automorphism of $A$ acting as the identity (resp., the inversion) map on $C_\infty$ (resp. on  $C_{p^\infty}$).  Then $G$ is {CF, but not BCF}.

\begin{thm}\cite{CDR}\label{CF-suff}
A  non-periodic {abelian-by-finite} group $G$ is {CF}  if and only if there is a normal abelian subgroup $A$ of finite index such that 
either $G$ acts by conjugation on $A$ as the power $\pm1$, or there is a $G$-series $1\le V\le A$ such that: 
\begin{enumerate}
\item $G/V$ is periodic and CF, 
\item $V$ is finitely generated free-abelian and $G$ acts by conjugation  on $V$
as the power $\pm1$.
\end{enumerate}
Moreover, $G$ is {BCF}  if and only if $G$ is CF and \\ 
\phantom{iiiiii}{\rm (3)}  there is $B\le A$ such that $G/B$ has finite exponent and $G$ acts by conjugation on $B$ as the power $\pm1$.
\end{thm}

%%%%%%%%%%%_____SUBSECTION___________

\subsection{Groups in which each subgroup is commensurable with a normal subgroup}\label{CN} 
In order to put the above results in a common framework, Casolo,  Dardano and Rinauro \cite{CDR} have recently considered the class of {\em CN-groups}, that is, groups in which each subgroup is commensurable with a normal subgroup. Clearly, a group is inertial if and only if 
(IN) below holds.  Moreover,  in view of Theorem \ref{BL}, a group is a  CN-group if and only if   (CN) below  holds.

\begin{definition}\label{def-CN}  Let $G$ be a group. Consider the follwing properties for $G$: \hfill  
\begin{itemize} %{enumerate}
  \item[(IN)]\ \ \  $\forall H\le G\ \forall g\in G\ \exists n\in\N\ \ |H:(H\cap H^g)|\le n$;
  \item[(CN)]\ \   $\, \, \, \forall H\le G\ \exists n\in\N\  \forall g\in G\ \ |H:(H\cap H^g)|\le n$; 
  \item[(BCN)]\  \  $\,	\,\, \exists n\in\N\ \forall H\le G\  \forall g\in G\ \ |H:(H\cap H^g)|\le n$.
\end{itemize} % {enumerate}
\end{definition} 

By Theorem \ref{BL}, {\em BCN-groups} (boundedly CN) are groups $G$ for which there is $n\in \N$ such that for each $H\le G$ there is $N\n G$ such that $|HN/H\cap N|\le n$. 

It is easy to show that the class of CN-groups is contained in the class of so-called {\em sbyf-groups}, that is, groups in which each subgroup $H$ contains a subnormal subgroup $S$ of $G$ such that the index $|H : S|$ is finite, say $H$ is {\em subnormal-by-finite}. Sbyf-groups have been studied by Heineken \cite{HH}  in 1996 and  Casolo \cite{sbyf} in 2016.

\begin{thm}\label{sbyf+} Let $G$ a locally finite sbyf-group. Then: 
\begin{enumerate} 
\item \cite{HH} $G$ is nilpotent-by-Chernikov; 
\item  \cite{sbyf} $G$ is (locally nilpotent)-by-finite.
\end{enumerate} 
\end{thm}

\begin{problem} Study sbyf-groups whose periodic quotients are locally finite.
\end{problem} 

 The next  result can be found in the forthcoming paper \cite{CDR}.

\begin{thm}\cite{CDR}\label{CN++}  Let $G$ be a CN-group. If every periodic quotient   of $G$ is locally finite, then $G$ is finite-by-abelian-by-finite.
\end{thm}

It is easy to see that for an abelian-by-finite group the properties CN and CF are equivalent. However, for every prime $p$ there is a nilpotent $p$-group with the property CN which is neither finite-by-abelian nor abelian-by-finite (\cite{CDR}).
\bigskip

\begin{thm}\cite{CDR}{ Let $G$ be a finite-by-abelian-by-finite group.
\begin{enumerate}
\item $G$ is \CN\ if and only if it is finite-by-\CF.
\item  $G$ is \BCN\ if and only if it is finite-by-\BCF.
\end{enumerate} }
\end{thm}

\vskip-1mm
It follows that if the group $G$ is periodic and finite-by-abelian-by-finite, then $G$ is BCN if and only if it is CN. 
 The more restrictive property BCN remains treatable when we consider the wider class of locally graded groups.  

\vskip-1mm
\begin{thm}\cite{CDR}  A locally graded BCN-group is finite-by-abelian-by-finite. \end{thm}

\vskip-1mm

The proof of  the above results will appear in \cite{CDR}. Here we just notice that in the proof of Theorem \ref{CN++} one  reduces first to the case when $G$ is locally nilpotent and soluble, by virtue of Theorem \ref{sbyf+}. Then one reduces to the case when $G$ is nilpotent, by using  Theorem \ref{Recall1} below and some technical lemmas. Then, finally, one uses techniques introduced by M{\" o}hres  \cite{M} for nilpotent groups.

\begin{definition}\label{def-AP} Let $A$ be an abelian group and $\Gamma\le Aut(A)$ a group of automorphisms of $A$. Consider the following properties for the pair $(\Gamma,A)$:  
\begin{itemize} 
\item[(\PP)]  $\forall H\le A\ \ H=H^\Gamma$; 
\item[(\AP)] $\forall H\le A\ \ |H/H_\Gamma|<\infty$; 
\item[(\BP)] $\forall H\le A\ \ |H^\Gamma/H|<\infty$; 
\item[(\CP)] $\forall H\le A\ \ \exists K=K^\Gamma\le A$ {such that $H \sim K$}, i.e. $H$, $K$ are  commensurable; 
\item[(\TP)] $\Gamma$ has $(\PP)$ on the factors of a series $1\le V\le D\le A$ whose members are $\Gamma$-invariant and 
  \begin{enumerate}
    \item   $V$ is free abelian of finite rank and $A/V$ is periodic, 
    \item  $D/V$ is divisible of finite rank and $\pi(D/V)$ is finite,
    \item the $p$-component of $A/D$ is bounded for each $p\in \pi(D/V)$ .  
  \end{enumerate}
\end{itemize}
\end{definition} 

\vskip-2mm
 
According to \cite{BLNSW}, $A$ is said to be  {\em $\Gamma$-hamiltonian}, when (\PP) holds. Obviously both (\AP)\ and (\BP)\ imply (\CP). Moreover,   (\AP), (\BP)\ and (\CP) {\em are equivalent,  provided that $A$ is abelian and $\Gamma$ is finitely generated, while they are in fact different in general even when both $A$ and $\Gamma$ are elementary abelian $p$-groups} (see \cite{DR1}).

\begin{thm}\label{Recall1}\  Let  $\Gamma$ be a group acting on an abelian group $A$. Then:\begin{enumerate} 
  \item  \cite{FGN} \ $\Gamma$ has $(\AP) $ on $A$ if and only if there is a $\Gamma$-invariant  subgroup $A_1$ such that $A/A_1$ is finite and $\Gamma$ has either $(\PP)$ or $($\TP$)$ on $A_1$;
  \item \cite{CasoloBP}    \ $\Gamma$ has $(\BP)$  on $A$ if and only if there is a $\Gamma$-invariant subgroup\ $A_0$\ such that $A_0 $ is\  finite and $\Gamma$ has either $(\PP)$ or $($\TP$)$  on $A/A_0$; 
  \item   \cite{CDR} $\Gamma$ has $(\CP)$ on $A$ if and only if  there are $\Gamma$-invariant subgroups $A_0\le A_1\le A$ such that $A_0$ and $A/A_1$ are finite and $\Gamma$ has either $(\PP)$ or $($\TP$)$  on $A_1/A_0$.
\end{enumerate} 
\end{thm}

%\vskip-1cm
\begin{Cor}\label{CorollarioCP} 
For a group $\Gamma$ acting on an abelian group $A$, the following are equivalent:
\begin{enumerate} 
  \item  $\Gamma$ has $(\AP)$ on $A/A_0$ for a finite  $\Gamma$-invariant subgroup $A_0$ of  $A$;
  \item   $\Gamma$ has $(\BP)$ on a  $\Gamma$-invariant subgroup $A_1$ of  $A$ of finite index;
  \item  $\Gamma$ has $(\CP)$ on $A$. 
\end{enumerate}
\end{Cor}

In closing this subsection we observe that the definition of inert subgroup suggests  considering also another natural subclass of the class of inertial groups.

\begin{problem}  Describe (soluble) groups which have the property:

$($BIN$)$\ \ \ $\forall g\in G \ \exists n\in\N\  \forall H\le G\ \ |H:(H\cap H^g)|\le n$.
\end{problem}

%%%%%%%%%%%%%%%______SUBSERTION____%%%%%%%%%%%%%%%%%%%------------------------

\subsection{Groups in which each subnormal subgroup is inert} 
The study of groups in which subnormal subgroups are subject to some restrictions is a standard problem in the theory of groups. 
The class of soluble groups in which subnormal subgroups are normal (usually called $T$-groups) and its generalizations has been studied  by many authors.  
 Gasch{\" u}tz \cite{G} treated  the finite case  in 1957. Then Robinson \cite{RT} treated the general case in 1964. Of course one has to consider classes of groups in which there are sufficiently many subnormal subgroups.

\begin{thm}\label{T}\cite{RT} Let $G$ be a soluble group in which subnormal subgroups are normal. Then $G$ is metabelian.  If  $G$ is finitely generated, then $G$ is finite or abelian.
\end{thm}

In 1985 Franciosi and de Giovanni considered  weaker  conditions.

\begin{thm}\label{IT}\cite{FG} Let $G$ be a soluble group in which subnormal subgroups are either normal or finite. Then the derived subgroup $G'$  of $G$ is nilpotent and $G''$ is cyclic of prime power order.
\end{thm} 

An application of this result can be found in \cite{D1}. On the other hand, in  1989 
Casolo \cite{CasoloBP}  studied the class of groups whose subnormal subgroups are nearly-normal. These are exactly those groups in which the relation of being nearly-normal is  transitive.  They have been called $T^*$-groups.

\begin{thm}\cite{CasoloBP} Let $G$ be a soluble  $T^*$-group. Then $G$ is finite-by-metabelian. If $G$ is finitely generated, then $G$ is abelian-by-finite.
\end{thm}

 A  stronger condition was considered by Franciosi and de Giovanni in 1993.

\begin{thm}\label{IT+}\cite{FG-LT} Let $G$ be a soluble infinite group in which subnormal subgroups are either normal or have finite index. Then $G$ is metabelian. \end{thm}

In 1995, Franciosi, de Giovanni and Newell considered a dual condition and studied $T_*$-groups, that is, groups whose subnormal subgroups are normal-by-finite. These are exactly those groups in which the relation of being  normal-by-finite is  transitive.

\begin{thm}\cite{FGN} Let $G$ be a soluble $T_*$-group. 
Then $G$ is metabelian-by-finite. Moreover, the derived subgroup $G'$ is abelian-by-finite. If $G$ is finitely generated, then $G$ is abelian-by-finite.
\end{thm}

In order to put these results in the same framework, we consider the class of  soluble group in which subnormal subgroups are (strongly) inert and report a result for finitely generated groups. Say that a group $G$ is {\em  semidihedral} on a torsion-free abelian subgroup $A$ if $G$ acts on $A$ by means of inertial automorphisms and $C_G(A) = A$.  This results to be equivalent to the fact that  $G=A\rtimes K$ where $A$ is  torsion-free abelian and $K$ acts faithfully on $A$ by means of inertial automorphisms (see \cite{DR5}). It is easy to verify that inertial automorphisms of a torsion-free abelian group are multiplications, according to Definition \ref{def-multiplication} and  Theorems \ref{M(A)} and \ref{caratterizzazione} below. Therefore, $K=K_0\times T$ where $K_0$ is free-abelian and $T$ has order at most $2$ (for details see  \cite{DR5}). The reader must be warned that the word {\em semidihedral} has been used also with a different meaning in other areas of group theory. 

\begin{thm} \cite{DR5}\label{sn-inert} Let $G$ be finitely generated soluble-by-finite group. Then each subnormal subgroup of $G$ is inert  if and only if $G$ has a finite normal subgroup $F$ such that $G/F$ is a semidihedral group.
\end{thm} 

Let us deduce a corollary generalizing the above stated fact that finitely generated soluble $T^*$- and $T_*$-groups are both abelian-by-finite.

\begin{Cor}
Let $G$ be a finitely generated soluble group. Then the following are equivalent:
\begin{itemize} 
  \item[1)] each subnormal subgroup of $G$ is commensurable with a normal subgroup;
  \item[2)]  $G$ has a finitely generated torsion-free abelian subgroup $A_0$ of finite index such that $G$ acts on $A_0$ by either the identity or the inversion map;
  \item[3)]  $G$ is CF, i.e., each subgroup is normal-by-finite.
\end{itemize}
\end{Cor}

\pf   Assume ($1$). Apply Theorem \ref{sn-inert} and   let $F$ be a finite normal subgroup  such that $G/F$ is semidihedral, thus $G/F=A/F\rtimes K/F$ where $A/F$ is  torsion-free abelian and $K/F$ acts faithfully on $A/F$ by means of inertial automorphisms (hence multiplications). Since each subgroup of $ A/F$ is commensurable with a subgroup
which is  normalized by $G$, these multiplications are either the identity or the inversion map. Since $K/F$ acts faithfully on $A/F$, we can conclude that $K/F$ and $G/A$ have  order $2$. Then $A$ is finitely generated. Moreover, $A$ is finite-by-abelian, hence FC (see \cite{Rm}). Then  $Z=Z(A)$ is a finitely generated abelian subgroup with finite index in $G$.
Now we can apply Theorem \ref{caratterizzazione} below (or Theorem 4 in \cite{DR1}) and deduce the existence of $A_0\le Z$ like in statement ($2$). 
%By Theorem \ref{CF-suff}, 
% ($2$) implies ($3$). Finally,  ($3$)  trivially implies ($1$).
Finally, it is clear that ($2$) implies ($3$) and ($3$) implies ($1$).
\QED

The next proposition shows that we are extending both  $T^*$ and $T_*$ in a natural way.

\begin{Prop}
For a group $G$ the following conditions are equivalent: 
\begin{enumerate} 
  \item  each subnormal subgroup of $G$ is commensurable with a normal subgroup;
  \item  the relation of being commensurable with a normal subgroup is transitive.
\end{enumerate} 
\end{Prop}

\pf Assume $(1)$ and let $H$ and $K$ be subgroups with the property that  there are subgroups $H_1$ and $K_1$ such that $H\sim  H_{1}\n K\sim K_{1}\n G$. We wish to prove that there is $H_2$ such that $H \sim H_2\n G$.  The subgroup $K_{2}:=(K\cap K_{1})_{KK_{1}}$has finite index in $KK_{1}$ and so $H_1\cap K_{2}$ has finite index in $H_1$ as $H_1\le KK_1$. Hence, $H\sim H_1\sim H_1\cap K_{2}$, 
according to Fact \ref{Fact}.
On the other hand $H_1\cap K_{2}\n K_2\n K_1 \n G$, hence $H_1\cap K_2$ is in turn commensurable with a normal subgroup. Thus, $(2)$ holds. The converse is obvious. 
\qed

\begin{problem}
Describe the structure of  soluble groups $G$ in which every subnormal subgroup is  commensurable with a normal subgroup.
\end{problem}

%%%%%%%%%%%%%%%__SECTION_3____%%%%%%%%%%%%%%%%%%%%%%%%%%%%%%%
\smallskip
\section{Groups of inertial automorphisms of a group}\label{Groups}

%%%%%%%%%

\subsection{The group of finitary  automorphisms of a group}  We start this section by  recalling a basic instance of inertial automorphisms. Let us extend the notation for the centralizer subgroup to denote the subgroup of  fixed-points of a set of  endomorphisms of a group.

\begin{definition}\label{auto-finitary}\ \  An endomorphism $\g$ of a group $G$ is called {\em finitary} if the subgroup\  $C_G(\g):=\{x\in G\ |\ x^\g=x\}$ has finite index in $G$. 
\end{definition}

Clearly, a group is an FC-group if and only if all its inner automorphisms are finitary. In any case, finitary automorphism 
 form a subgroup, say $\FAut(G)$, of the group of all automorphisms of $G$.  The group  $\FAut(G)$ seems to have been studied for the first time in 1975  in \cite{Zal} by A. E. Zalesskii (see Theorem \ref{BeShv}(2) below) who  called {\it bounded automorphisms} what we here call  finitary automorphisms (following  \cite{W}) .

The case when the underlying group is abelian, has been considered independently by   Wehrfritz in 2002.

\begin{thm}\label{W1} \cite{W} If $A$ is an abelian group, then the group  $\FAut(A)$  of all finitary automorphisms of $A$  is locally finite.
\end{thm}

 Imposing some  additional restraints on the underlying abelian group $A$,  one can describe the group $\FAut(A)$. For example,  
observe that if $A$ is torsion-free (or divisible), then $\FAut(A)=1$, since for each  $\g\in\FAut(A)$ the subgroup  $A/ker(\g-id)\iso \im(\g-id)$ is 
simultaneously finite and torsion-free, hence trivial.

Now we report two results which suggest a construction for $\FAut(A)$ for some abelian groups $A$. If $G$ is any group and   $N$ a normal subgroup of $G$, denote by $\St(G,N)$ the stability group in $\Aut(G)$ of the series $G\ge N\ge 1$, that is, the set of automorphisms $\g\in \Aut(G)$ such that $N$ is a $\la \g\ra$-invariant normal subgroup of $G$ and $\g$ acts as the identity map on both $N$ and $G/N$. Then it is well-known that  $\St(G,N)$ is abelian and, if $G$ is abelian,  
we have  $\St(G,N)\iso \Hom(G/N,N)$.

\begin{thm}\label{FAut(critico)}  \cite{DR3} Let $A$ be an abelian $p$-group such that $D:=D(A)$ has finite rank and $A/D$ is bounded. 
Then,  $\Sigma:=\St(A,D)$ is a bounded abelian $p$-group and there is a subgroup $\Phi\iso \FAut(A/D)$ of    $\FAut(A)$ such that 
$\FAut(A)=\Sigma\rtimes \Phi$ 
 where the automorphisms induced by $\Phi$ via conjugation on $\Sigma$ are finitary and this action is faithful.\end{thm}

\begin{thm}\label{IAut(Asplits)} \cite{DR3}
 Let $A$ be an abelian group and $T=t(A)$. If 
$A/T$ is  finitely generated (resp., $T$ is bounded), then $\Sigma:=\St(A,T)$ is a periodic (resp., bounded) abelian group and there is a subgroup  $\Phi_1\iso \FAut(T)$ of   $\FAut(A)$ such that $\FAut(A)=\Sigma\rtimes \Phi_1$, 
 where  $\Phi_1$ induces  finitary automorphisms by conjugation on $\Sigma$ .

In the case when $A/T\ne 0$ is finitely generated, the above action by conjugation is faithful. \end{thm}

Concerning the general case in which the underlying group  need not be abelian the following holds. 

\begin{thm}\label{BeShv}\cite{BS1,BS2, MR, Zal} Let $G$ be any group and $\FAut(G)$ be the group of all finitary automorphisms of $G$. Then:
\begin{enumerate} 
 \item $\FAut(G)$  is  abelian-by-(locally finite);
  \item   $\FAut(G)$  is locally (centre-by-finite), hence (locally finite)-by-abelian; 
  \item   $\FAut(G)$ is locally finite if and only if $\FAut(G)\cap Inn(G)$ is locally finite.
\end{enumerate}
\end{thm}

 Statement (1) has been shown in 2009 by Belyaev and Shved in \cite{BS1,BS2}. Part (2) follows from  Zalesskii's investigation in \cite{Zal} and  
implies that the elements of $\FAut(G)$ of finite order form a locally finite subgroup. Finally, part (3) follows from results in \cite{MR} by  Menegazzo and Robinson in  1987.

Recently,  Shved has obtained the following result. Recall that a finitary linear group is a group of automorphisms of a vector space which fix poitwise a subspace of finite codimension.

\begin{thm}\label{Shv}\cite{Shved} Let $G$ be any group, then $\FAut(G)$ has a normal  series 
$$1 \le \Gamma_1 \le \Gamma_2\le \FAut(G)$$  such that:
\begin{enumerate}
  \item  $\Gamma_1$  is nilpotent of class al most 4;
  \item   $\Gamma_2/\Gamma_1$ is periodic and has an ascending central series of type at most $\omega$; 
  \item  $\FAut(G)/\Gamma_2$ is embeddable into a direct product of finitary linear groups over prime fields.
\end{enumerate}
\end{thm}
\noindent In \cite{Shved} more details are given in the case when $G$ is a $p$-group (finitary automorphisms are called {\it virtually trivial}  in that paper).

\bigskip 

Let us  consider now a particular type of finitary automorphism that we will use in \S \ref{strongly-inertial+}.

\begin{definition}\label{stronglyfinitary}  An automorphism $\g$ of a group $G$ is called {\em strongly finitary} if  there is a finite $\la \g\ra$-invariant normal 
subgroup $N$ of $G$ such that $\g$ acts trivially on $G/N$. 
\end{definition}

In the above notation, since the correspondence $x \mapsto[x,\g]$ gives an injection of a right transversal of  $C_G(\g)$ in $G$ into $N$,  it results $|G:C_G(\g)|\le |N|$. Thus, the strongly finitary automorphisms of $G$ form a subset  $s\FAut(G)$ of $\FAut(G)$. Moreover, if $G$ is abelian we have $s\FAut(G)=\FAut(G)$. In fact, if  
$\g\in\FAut(G)$ then  $A/ker(\g-id)\iso \im(\g-id)=:N$ is finite and $\g$ acts on $G/N$  as the identity map. On the other hand, in the dihedral group $\la x,y \ |\  x^2=1, y^x=y^{-1}\ra$ the inner automorphism induced by $y$ is finitary but not strongly finitary. Thus we can have $s\FAut(G)\ne\FAut(G)$. 

Then let us state an easy extension of  Theorem \ref{W1}.

\begin{thm}\label{W+}  Let  $G$ be a group. Then the group  $s\FAut(G)$  of all strongly finitary automorphisms of $G$  is locally finite and coincides with the subgroup of $\FAut(G)$ formed by  the elements of finite order.
\end{thm}

\pf Let $\g_1,\dots,\g_n\in s\FAut(G)$. If $\g_i$ acts trivially on $G/N_i$, then  the subgroup $N:=N_1\cdot\ldots \cdot N_n$ is finite, normal in $G$ and $\la \g_i\ra$-invariant for each $i$. Then each $\g_i$ acts trivially on $G/N$ and $s\FAut(G)$ is a subgroup of $\FAut(G)$. Moreover,  if $\Gamma=\la\g_1,\dots,\g_n\ra$,  then the subgroup   $C:=C_\Gamma(N)$ has finite index in $\Gamma$ and stabilizes the  series  $1\le N \le G$. Thus $C$ is abelian and finitely generated. 
 With $m= |N|$, for each $\g_i$ and $a\in G$ we have $[\g_i^m,a]=[\g_i,a]^m=1$. Then $C$ and $\Gamma$ are finite. Hence $s\FAut(G)$ is locally finite.

On the other hand, if $\g\in\FAut(G)$ and $\Gamma_1:=\la \g\ra$ is finite, consider 
$G\rtimes \Gamma_1$ as a subgroup of the holomorph group of $G$. By Dicman's Lemma (see \cite{Rm}), $\Gamma_1^{G}$ is finite. Let $N:=[G,\Gamma_1]\le \Gamma_1^{G}$, then $N$ is a  finite normal subgroup of $G$ and $\Gamma_1$ acts trivially on $G/N$. Hence $\g\in s\FAut(G)$.\qed

%%%%%%%%%%%%%%%%%%%%%%%%%%%%%%

\subsection{The group of automorphisms of a group fixing all but finitely many subgroups} 
Power automorphisms are trivially inertial and have been studied by many authors, mainly by C.D.H. Cooper \cite{Coo} in 1968. 

\begin{thm}\cite{Coo} The set $\PAut G$ of all power automorphisms of a group $G$ is a normal, abelian, residually finite subgroup of the full automorphism group $\Aut(G)$ of $G$.
\end{thm}

The subgroup  $N(G)$ of all elements of a group $G$ which induce by conjugation a power automorphism is called the {\em norm}, it coincides with  the intersection of all normalizers of subgroups of $G$. A group $G$ is a Dedekind group if and only if $G=N(G)$. Clearly we have $Z(G)\le N(G)$. Furthermore, by a classical results by Wielandt, we have $N(G)\le Z_2(G)$ where $Z_2(G)/Z(G):=Z(G/Z(G))$ is the so-called second centre of $G$. This follows also from the fact 
that power automorphisms are central,  that is, they act trivially on the central factor group $G/Z(G)$ (see \cite{Coo}). Equivalently, an automorphism of $G$ is central if and only if it centralizes the group $Inn(G)$ of inner automorphisms (in the full automorphism group $\Aut G$ of $G$) .

In 1993 Cutolo  \cite{Cut}   studied  the group $\QAut(G)$ of  the  so-called {\em quasi-power} automorphisms of a group $G$, that is,  automorphisms which fix all but finitely many subgroups of $G$. 

\begin{thm}\cite{Cut}  Let $G$ be an infinite group. Then $\QAut(G)$ is  abelian and residually finite.  

If  $\QAut(G) \ne  \PAut (G)$, then the subgroups of $G$ which are not fixed under the action of $\QAut G$ generate a finite characteristic subgroup $F$ such that $G/F$ is a finite extension of a $p$-subgroup. In particular, $G$ is periodic. Moreover, every infinite subgroup of $G$ which is locally finite is Pr\" ufer-by-finite.
\end{thm}

\subsection{The group of automorphisms of a group fixing all infinite subgroups} 
If $IN(G)$ is the subgroup of the group $G$ defined as  the intersection of all normalizers of infinite subgroups of $G$ we have clearly $Z(G)\le N(G)\le IN(G)$. In 2000,    Beidelman and Heineken in \cite{BH} give many results about this subgroup. We report just a couple of facts.

\begin{thm}\cite{BH} 
If $G$ is a non-periodic group, then $IN(G)$ is abelian. Moreover, if $IN(G)$ is non-periodic, then $IN(G)=Z(G)$.
\end{thm}

In 1990 Curzio, Franciosi and de Giovanni   \cite{CFG} had considered the group I-$\Aut(G)$ of all automorphisms $\g$ of $G$ such that $H^\g=H$ for each infinite subgroup $H$ of $G$. These automorphisms are clearly inertial  and I-$\Aut(G)$ is a normal subgroup of $\Aut(G)$ containing the groups $\PAut(G)$ and $\QAut(G)$, see \cite{Cut}. We summarize their results as follows.

\begin{thm}\cite{CFG} Let $G$ be a group.\begin{enumerate}
  \item If $G$ is non-periodic, then  {\rm I}-$\Aut(G)$ is abelian.
  \item There are infinite hypercentral $p$-groups $G$ for which {\rm I}-$\Aut(G)$ is not abelian.
  \item If  $G$ is   locally finite with no infinite section which are a simple groups and  {\rm I}-$\Aut(G)$ is not abelian, then $G$ is Chernikov.
  \item If $G$ is (locally radical)-by-(finite) and  {\rm I}-$\Aut(G)\ne \PAut(G)$, then $G$ is Chernikov.
  \item If $G$ is either FC or nilpotent and  {\rm I}-$\Aut(G)\ne \PAut(G)$, then $G$ is Pr{\"u}fer-by-finite.
\end{enumerate}
\end{thm}

Recall that a {\em radical} group is a group which has an ascending series with locally nilpotent factors.

In \cite{BH} Beidelman and Heineken have continued the study of the group I-$\Aut(G)$ and have given further information. For example, they show that in most cases the equality I-$\Aut(G) = \PAut(G)$ actually holds and give a complete characterization of Chernikov groups satisfying I-$\Aut(G) \not= \PAut(G)$. Their results are rather technical. For the sake of brevity, we recall only the following one and we refer the reader to \cite{BH} for the rest. 

\begin{thm}\cite{BH} 
If $G$ is an infinite locally finite group, then {\rm I}-$\Aut(G)$ is metabelian.
\end{thm}

Recently, 
%concerning non-periodic groups, 
De Falco, de Giovanni, Musella and Sysak in \cite{FGMS17} have introduced the concept of {\it weakly power automorphism}, that is an 
 automorphism $\gamma$ of a group $G$ such that $H^\gamma=H$ for each non-periodic subgroup $H$ of $G$. Denote by $\WAut(G)$ the group of all such automorphisms. Clearly  $\WAut(G)$ contains I-$\Aut(G)$.

\begin{thm} \cite{FGMS17} Let $G$ be a non-periodic  group. 
\begin{enumerate}
 \item The group $\WAut(G)$ is abelian.
 \item If $G$ is generated by non-periodic elements, then $\WAut(G)=\PAut(G)$.
 \item If $G'$ is metabelian-by-finite, then $\WAut(G)=\PAut(G)$.
 \item There exists a soluble non-periodic group (with derived length $4$) such that  $\WAut(G)\ne \PAut(G)$.
\end{enumerate}
\end{thm}

%%%%%%%%%%%%%%%%%%%%%%%%%%%%%%

\subsection{The group of  inertial automorphisms of an abelian group}\label{IAut(A)} 
As we observed above, the  inverse of an inertial automorphism need not to be inertial.
On the other hand, if $G$ is any group and $\g\in Aut(G)$, then both $\g$ and $\g^{-1}$ are inertial if and only if $H\sim H^\g$ for each  $H\le G$. Therefore, we introduce a couple of definitions.

\begin{definition}\label{bi-inertial}  If $\g$ is an automorphism of a group $G$ we say that:
\begin{enumerate}
  \item   $\g$ is  {\em bi-inertial} if  both $\g$ and $\g^{-1}$ are inertial;
  \item  $\g$ is  {\em almost-power} if every subgroup $H$ of $G$ contains a subgroup $K$ of finite index such that $K^\g=K$.
\end{enumerate}
\end{definition} 

 According to the terminology introduced in Section \ref{normal-by-finite}, an automorphism $\g$ is almost-power if and only if the group $\la \g\ra$ has (\AP). Clearly an almost-power automorphism is bi-inertial. Moreover,  a periodic inertial automorphism is almost-power. The definition of almost-power automorphism was introduced in  \cite{FGN} where it was shown that the set of almost-power automorphisms of an abelian group $A$ is a subgroup of $\Aut(A)$. By Fact \ref{Fact}, it follows that {\em the set of all bi-inertial automorphisms of any group $G$ is a subgroup of $\Aut(G)$}  which contains  all almost-power automorphisms   and seems to be an object of interest.  However, in the general case, this group may contain even a Tarski  group, hence it can be  too complicated. Also a centreless FC-group may be regarded as a group of bi-inertial automorphisms. On the other hand, if the underlying group is abelian the picture is clearer, as next theorem shows. 
As we announced, we will focus attention on 
the group $\IAut(A)$ {\em generated by the inertial automorphisms} of an abelian group $A$ and report results from \cite{DR3}.  We warn that in  \cite{DR1} by inertial one means what here is called bi-inertial. The terminology of the present paper is consistent with  \cite{DR2,DR3}.

\begin{thm}\label{autoinerziali}\cite{DR1,DR2} Let $\g$,\ $\g_1$ be  automorphisms of an abelian group $A$.
\begin{enumerate}
\item If  $\g,\g_1$ are both inertial, then the commutator $[\g,\g_1]$ is finitary.
\item The group $\IAut(A)$ coincides with the set of the products  $\g\g_1^{-1}$ where both $\g$,$\g_1$ are inertial. 
\item If $r_0(A)<\infty$, then $\g$ is inertial if and only if $\g$ is bi-inertial.
\item If $r_0(A)=\infty$, then $\g$ is bi-inertial if and only if $\g$ is almost-power.
\end{enumerate}
 \end{thm}

Recall that $\IAut(A)$ contains the group $\FAut(A)$, which is locally finite, by Theorem \ref{W1}.  In the case $A$ is periodic, we can give many details on the group  $\IAut(A)$ in next theorem. If $\pi$ is a set of primes, we  regard $\FAut(A_\pi)$ as naturally embedded in $\FAut(A)$. 

\begin{thm}\cite{DR3}
Let $A$ be a periodic abelian group, then $\IAut(A)=\PAut(A)\cdot \Delta\cdot \FAut(A)$
where $\Delta$ is a direct product of cyclic  groups  (with  $\Delta=1$, 
if $A$ is reduced). Moreover, all inertial automorphisms of $A$ are almost-power.

Therefore, $\IAut(A)$ is central-by-(locally finite), hence (locally finite)-by-abelian.

Furthermore, there is a set $\pi$  of primes such that
$\Delta\cdot \FAut(A)=\FAut(A_\pi)\times (\Sigma \ltimes\ \mathcal{I})$, 
 where $\Sigma$ is an abelian $\pi'$-group with bounded primary components, $\Delta\le\mathcal{I}$ and
 $\mathcal{I}$ acts faithfully on $\Sigma$ by inertial automorphisms.
\end{thm}

To treat the  general case, we list results from \cite{DR3} in a unique theorem. 
Let $\IAut_1(A)$ be the group of bi-inertial automorphisms
centralizing $A/t(A)$.

\begin{thm}\cite{DR3} Let $A$ be an abelian group.\begin{enumerate}
  \item   
 If $A$ is torsion-free, then $\IAut(A)$ is abelian.
  \item 
 $ \FAut(A)\le \IAut_1(A)$ and equality holds if
$r_0(A)=\infty$.
  \item  
 $\IAut(A)=\IAut_1(A)\times Q(A)$ where $Q(A)$ is isomorphic to a multiplicative group of  rational numbers.
  \item   $\IAut_1(A)\times \{\pm 1\}$ coincides with the group of almost-power automorphisms of $A$.
  \item  $\IAut(A)$ is (locally finite)-by-abelian.
  \item   There is a normal subgroup $\Gamma\le \IAut_1(A)$ such that:\\
 i)\ \ $\IAut_1(A)/\Gamma$ is locally finite;\\
ii)\ $\Gamma$ acts by means of power automorphisms on its derived subgroup $\Gamma'$, which is a periodic abelian group.  
  \item    $\IAut(A)$ is metabelian-by-(locally finite).
  \item   If $A=\Z(p^\infty)\oplus \Z$, then $\IAut(A)$ is not (locally nilpotent)-by-(locally finite).
\end{enumerate}
\end{thm}

 Recall that a metabelian group  $G$ in which each subgroup of the derived subgroup  $G'$ is normal in $G$ is  said to be a {\em KI-group}. Such groups have been studied by Subbotin in a series of papers (see \cite{S-ZD} and the references therein).
 A further instance of KI-group is given in the following theorem. By $Fit(G)$ we denote the subgroup generated by all normal nilpotent subgroups of $G$

\begin{thm}\cite{DF1} Let $G$ be a group and $\Aut_{sn}(G)$ be the group  of all automorphisms $\g$ of $G$ such that $H^\g=H$  for each subnormal subgroup $H$ of $G$. If $G$ is soluble and  either $G$ is nilpotent-by-(finitely generated) or $\Fit(G)$ is non-periodic, then  $\Aut_{sn}(G)$ is a KI-group.
\end{thm}

 Let us recall that, to show that a group of automorphisms is a KI-group one can use arguments like the following lemma which is inspired by the celebrated argument introduced by P. Hall to prove that the stability group of a finite series must be nilpotent (see \cite{CP,DF2}).
 
\begin{lemma}\label{LemmaDF2} \cite{DF2} Let $A$ be an abelian group, $\sigma,\g\in \Aut(A)$ and $m_{1}, m_{2}\in\Z$. If $\sigma$ stabilizes a series $0\le A_1\le A$ and   $a_1^\g=a_1^{m_{1}}$ for each $a_1\in A_{1}$, $[ a ]^{\g^{-1}}=[ a] ^{ m_{2}}$ for each $[a]\in A/A_{1}$, then $\sigma^\g=\sigma^{m_1m_2}$.
\end{lemma}

One can give a detailed description of $\IAut(A)$ in some cases. 

\begin{thm}\label{IAut(AsplitsTbounded)} \cite{DR3}
 Let $A$ be an abelian group with $r_0(A)<\infty$  and  $T:=t(A)$ be bounded. Then there are an abelian subgroup $\Sigma:=\St(A,T)$ of finite exponent and  a subgroup $\Gamma_1$ of $\IAut_{1}(A)$    such that $\IAut_{1}(A)=\Sigma\rtimes\Gamma_1$, 
where $\Gamma_1\iso \IAut(T)$  induces via conjugation on $\Sigma$ inertial automorphisms.
\end{thm}

\begin{thm}\label{IAut()}\cite{DR3} Let $A$ be a non periodic abelian group  and  $T:=t(A)$ 
with $A/T$  finitely generated. Then there are a periodic abelian  group $\Sigma:=\St(A,T)$ and a subgroup $\Gamma_1$   of \ $\IAut_{1}(A)$  such that
 $\IAut_{1}(A)=\Sigma\rtimes\Gamma_1$, where    $\Gamma_1\iso \IAut(T)$ induces via conjugation on $\Sigma$ inertial automorphisms and the action is faithful. 

If in addition $T$ is not bounded, then $\IAut_{1}(A)$ is not nilpotent-by-(locally finite). Further, if $A_{2'}$ is unbounded, then $\IAut_{1}(A)$ is not even (locally nilpotent)-by-(locally finite).
\end{thm}

%%%%%%%%%%%%%%%%%%%%%%%%%

\subsection{The group of strongly inertial automorphisms}\label{strongly-inertial+}
As observed above, groups of inertial automorphisms can have a complicated structure. This suggests to consider a stronger condition.

\begin{definition}
An automorphism $\g$ of a group $G$ is called {\em strongly inertial} if and only if:
$$
\forall H\le G\ \ |\langle H,H^\g\rangle : (H\cap H^\g)|<\infty.
$$
\end{definition}

This is equivalent to the fact that both $H$ and $H^\g$ have finite index in their join $\langle H,H^\g\rangle$.  Clearly a strongly inertial automorphism is inertial. Moreover, we have  $  |\langle H,H^\g\rangle : (H\cap H^\g)|=|\langle H^{\g^{-1}},H\rangle: (H^{\g^{-1}}\cap H)|$. Therefore, if $\g$ is strongly inertial, then $\g^{-1}$ is strongly inertial as well. Then a strongly inertial automorphism is bi-inertial. Moreover, when the underlying group is abelian, an inertial automorphism is strongly inertial and the two notions coincide.

\begin{proposition} Strongly inertial automorphisms of any group $G$ form a subgroup $s\IAut(G)$ of $\Aut(G)$ containing the subgroup $s\FAut(G)$ of strongly finitary automorphisms of $G$.
\end{proposition}

\pf If $\g_1$ and $\g_2$ are strongly inertial,   then the subgroup $H\cap H^{\g_1\g_2}$ has finite index in $H$ for each $H\le G$, since $\g_1\g_2$ is inertial, as we said above. Moreover, $H$ has finite index in $H_1:=\langle H,H^{\g_1}\rangle$ which in turn has finite index in $\langle H_1,H_1^{\g_2}\rangle$ which contains  $\langle H,H^{\g_1\g_2}\rangle$. For what concerns $s\FAut(G)$ see Theorem \ref{W+} and note that, if $\g\in s\FAut(G)$ acts trivially on $G/N$, then $\la H,H^\g\ra\le HN$.
\qed
\medskip

 In the dihedral group $\la x,y \ |\  x^2=1, y^x=y^{-1}\ra$ the inner automorphism induced by $y$ is finitary hence inertial, but not strongly inertial as $\la x, x^y\ra$ is infinite.

Clearly, {\em a group is strongly inertial if and only if each inner automorphism is strongly inertial}.

\begin{problem} Study the group $s\IAut(G)$, when $G$ is any group. \end{problem}

%%%%%%%%%%%%%%______SECTION_4___________%%%%%%%%%%%%%%%

\medskip 
\section{Endomorphisms of an abelian group}\label{abelian}

\subsection{Inertial endomorphisms of an abelian group}\label{endoabelian} 
From \S \ref{basic} recall that an endomorphism $\p$  of an abelian group $A$ such that $H+\p(H)/H$ is finite  for each subgroup $H$ of $A$ is called an inertial endomorphism. 
Recall also that an endomorphism $\p$ of $A$ is called finitary if the subgroup $\{x\in A\ |\ \p(x)=x\}$ has finite index in $A$.

Note that finitary endomorphisms are trivially inertial. More generally, if an endomorphism $\p$ of an abelian  group $A$  acts as an inertial endomorphism  modulo a finite subgroup of $A$ (or on a subgroup of finite index), then $\p$ is inertial on the whole $A$. 
To consider other instances of inertial endomorphisms we need to introduce some terminology.

\begin{definition}\label{def-multiplication}
An endomorphism $\p$ of an abelian group $A$ is called a  \emph{multiplication} if one of the following holds:
\begin{enumerate}
\item  $A$ is periodic and  $\p(H)\subseteq H$ for each subgroup $H$ of $A$; 
\item  $A$ is not periodic and there are ${m},{n}\in \Z$ such that $A=nA$, 
$A_{\pi(n)}=0$ 
and $\p(nx)=mx$ for all $x\in A$
\end{enumerate}
\end{definition}

If (1) holds, as we noticed in \S \ref{TIN}, for a power automorphism $\g$,  there exist a $p$-adic integer $\alpha$ such that we have $\g(a)=\alpha(a)$ for all $a\in A_p$ (here we are in additive notation for $A$ and $\alpha(a)$ stands for $na$, where $n$ is any integer congruent to $\alpha$ modulo the order of $a$). We can say 
that $\p$ acts as the multiplication by $\alpha$ on $A_p$. In case (2)  the abelian group $A$ has a natural structure of $\Q^{\pi(n)}$-module and we say that $\p$ acts as the multiplication by $\frac{m}{n}\in\Q$.  By abuse  of notation, we write  $\p=\frac{m}{n}$. Note that we use the word ``multiplication" in a way different from \cite{F}. Ours are in fact ``scalar" multiplications. 

\begin{thm}\cite{DR2,DR4}\label{M(A)} Let $A$ be an abelian group.
\begin{enumerate}
\item  The multiplications of  $A$ form a ring $M(A)$ and commute with any endomorphism. 
\item  If $A$ is torsion-free, then all inertial endomorphisms of $A$ are multiplications.
\item If $r_0(A)<\infty$, then all multiplications are inertial.
\end{enumerate}
\end{thm}

In order to study  inertial endomorphisms of torsion abelian groups, note first  that 
a standard argument shows  that {\em an endomorphism $\p$ of an abelian torsion group $A$  is inertial if and only if $\p$ is inertial 
on all primary components $A_p$ and multiplication on all but finitely many of them}.
So that we are reduced to consider inertial endomorphisms of abelian $p$-groups. 

\begin{thm}\label{TEOREMA_PERIODICO+}\cite{DR2} 
Let  $\p$ be an endomorphism of an abelian $p$-group $A$. Then $\p$ is inertial if and only if there is a $\p$-invariant subgroup $A_0$ of $A$ of finite index such that either (1) or (2) holds:
\begin{itemize}
 \item[(1)] $\p$ acts as a multiplication on $A_0$;
\item[(2)] $A_0=B\oplus D$ where $B$ is bounded, $D$ is divisible  of  finite rank,
$B$ and $D$ are
 $\p$-invariant and  $\p$ acts as a  multiplication (by possibly different $p$-adics)  on both $B$ and $D$.
\end{itemize} 
If $\p$ is inertial, then  there is $m\in\N$ such that for each $X\le A$ there are  $\p$-invariant  subgroups $X_*,X^*$ of $A$   such that $X_*\le X\le X^*$ and $|X^*/X_*|\le m.$
\end{thm}

For the case of non-periodic groups we have the following.

\begin{thm}\label{PropMixENDO}\cite{DR2}%caratterizzazione
Let  $\p$ be an endomorphism of a  non-periodic abelian group $A$. Then $\p$ is inertial if and only 
 if either (1) or (2) holds:\begin{itemize}
 \item[(1)] 
 there is a $\p$-invariant  subgroup $A_0$ of $A$ of finite index such that
$\p=m\in\Z$ on $A_0$;
 \item[(2)] there are $\frac{m}{n}\in\Q$, $r\in\N$ and a  $\p$-invariant subgroup  $V\iso \Q^{\pi(n)}\oplus\ldots\plus \Q^{\pi(n)}$ ($r$ times) of $A$ such that 
\begin{itemize}
  \item[i)]  $\p$ acts as the multiplication by $\frac{m}{n}$ on $V$; 
  \item[ii)]  the factor group $A/V$ is torsion and $\p$ induces an inertial 
    endomorphism on $A/V$;
 \item[iii)] $A_{\pi(n)}$ is bounded.
\end{itemize} 
\end{itemize} 
\end{thm}

We give also a more detailed (and technical) characterization for the action of finitely many  inertial endomorphisms.

\begin{thm}\label{caratterizzazione}\cite{DR2}.  Let $\p_1,\dots,\p_l$ be finitely many endomorphisms of an abelian  group $A$. Then each
$\p_i$ is inertial if and only if there is a subgroup $A_{0}$ of $A$ of  finite index which is $\p_i$-invariant for each $i$ and such that 
either (1) or (2) holds:\begin{itemize}
 \item[(1)]  $\p_i$ acts as multiplication on $A_{0}$ by some $m_i\in\Z$ for each $i$;
 \item[(2)]  $A_{0}=B\oplus D\oplus C$ where $B,C,D$ are $\p_i$-invariant for each $i$ and there exist finite sets of primes $\pi$, $\pi_1$ such that $\pi\subseteq \pi_1$  and
\begin{itemize}
\item[i)]\ $B\oplus D$ is the $\pi_1$-component of $A_{0}$, where $B$ is bounded and $D$ is a divisible $\pi'$-group of finite rank;
\item[ii)] $C$ is a $\pi$-divisible subgroup with $C_{\pi}=\{0\}$  and there is a subgroup $V\iso \Q^{\pi}\oplus\dots\oplus\Q^{\pi}$ (finitely many times) which is  $\p_i$-invariant for each $i$ and such that $C/V$ is a 
$\pi_1$-divisible $\pi'$-group; 
\item[iii)]each $\p_i$ acts by  multiplication  on each $B$, $D$, $V$, $C/V$ where $\p_i$ acts as $m_i/n_i\in\Q$ on $V$ and on all $p$-components of $D$ such that the $p$-component of $C/V$ is infinite and $\pi=\pi(n_1\cdots n_l)$.
\end{itemize}
\end{itemize}
\end{thm}

%%%%%%%%%%%%

\subsection{The ring of  inertial endomorphisms of an abelian group}\label{Ring} 

We have already seen that inertial endomorphisms of any abelian group $A$ form a ring, denoted by $\IEnd(A)$, containing the ideal $F(A)$ of  endomorphisms  with finite image.

In \cite{DR4}  one finds a rather satisfactory description of the whole ring  $\IEnd(A)$. To report it in the general case we would need to introduce many  tools which are too technical for the present paper. Here we recall the perhaps most important results  and refer the reader to \cite{DR4}. 

\begin{thm}\cite{DR2} 
The ring $\IEnd(A)$ of inertial endomorphisms of any abelian group $A$ is commutative modulo its ideal $F(A)$ formed by the endomorphisms of $A$ with finite image.
\end{thm}

In the case $r_0(A)<\infty$, one may wonder if the ring $\IEnd(A)$ is generated by the subrings $M(A)$ and $F(A)$. Actually, this is not the case, even for $p$-groups. However for $p$-groups we have the following theorem. Say that an abelian $p$-group is \emph{critical} when the maximum divisible subgroup $D(A)$ has positive finite rank and $A/D(A)$ is bounded but infinite.

\begin{thm}\label{ENDO_mod_F} \cite{DR4}  Let $A$ be an abelian $p$-group.\begin{enumerate}
\item If $A$ is non-critical then $\IEnd(A)= F(A)+M(A)$.
\item If $A$ is critical, then $\IEnd(A)= M(A)+F(A)+N$, where $N\iso\Z(p^e)$ and $p^{e}$ is the exponent of $A/D(A)$.\end{enumerate}
\end{thm}

%%%%%%%%%%%%

\subsection{Inertiality for linear mappings}\label{linear}  Recall that in 1988 Hall and Phillips in \cite{JH,Ph} considered  finitary linear maps (i.e. endomorphisms) of the vector space $V$ to be those (invertible) linear maps $\beta$ such that the  subspace $(\beta-id)V$ is finite-dimensional. This is clearly equivalent to the fact that $\beta$ acts as the identity map on a subspace of finite codimension. The study of groups of finitary linear maps has received much attention in the literature (see \cite{KS,MPP,W+}, and the bibliography therein for example). Along the same lines, Dardano and Rinauro \cite{DR4} have considered the corresponding inertial properties for linear maps. 

\begin{thm}\label{TEOREMA_LINEARE--}\cite{DR4} Let $\beta$ be an endomorphism of an infinite dimensional  vector space $V$. \begin{enumerate}
\item Each subspace $H$ of $V$ has finite codimension in $H+\beta(H)$  if and only if $\beta$ acts as a scalar multiplication on a finite codimension subspace of $V$.
\item For each subspace $H$ of $V$ the codimension of $H \cap \beta(H)$ in $H$ is finite  if and only if $\beta$ acts as multiplication by a non-zero scalar on a finite codimension subspace.
\end{enumerate}
\end{thm}

 \begin{Cor} The ring of endomorphisms $\beta$ of a vector space $V$ such that each subspace $H$ has finite codimension in $H+\p(H)$ is the sum of the subring of scalar multiplications and the subring of finitary endomorphisms.
  \end{Cor}

%%%%%%%%%%%%%%%%%%%%%%%%%%%%%%%%%%%%%%

%______INIZIO SEZIONE 5__FULLY____________________

%\medskip
\section{Fully inert vs fully invariant subgroups of abelian groups}\label{Fully}%%%%%%%%

%%%%%%%%%%%%%%%%%%%%
\subsection{The main problem and a conjecture}\label{reticoli}%%%%%%%%%%%

The main object of this section is the study of solutions to the general Problem \ref{3reticoli} on the interrelations between the following sublattices (of the lattice of subgroups of an abelian group $A$) that we introduced at the end of \S\ref{uniformly}, see Proposition \ref{uniformly}:
$$
\mathcal Inv(A)\ \subseteq \ \mathcal Inv^\sim(A)\ \subseteq \ \mathcal I_u(A)\ \subseteq \ \mathcal I(A). \eqno(*)
$$ 

Two natural more specific problems arise right away: 

\begin{problem}\label{Problem:I} {$ $ }
\begin{itemize}
   \item[(1)] Describe $\mathcal Inv^\sim(A)$ for an abelian group $A$, i.e., determine when is a fully inert subgroup of $A$ commensurable with a fully invariant subgroup of $A$.
   \item[(2)] Characterize the class $\mathfrak I$ of abelian groups $A$ with $\mathcal I(A) = \mathcal Inv^\sim(A)$, i.e., such that every fully inert subgroup of $A$  is commensurable with some fully invariant subgroup of $A$. 
\end{itemize}
\end{problem}

Let us involve now also the lattice $\mathcal I_u(A)$ in the chain (*). This gives rise to other two classes:  
the class $\mathfrak I_*$ (resp.,  $\mathfrak I_u$) of all groups $A$ where the equality $\mathcal Inv^\sim(A)= \mathcal I_u(A)$ 
(resp.,  $\mathcal I_u(A) =  \mathcal I(A)$) holds. Clearly, $\mathfrak I = \mathfrak I_*\cap \mathfrak I_u$, as $\mathfrak I$  is the class of all abelian groups $A$  such that $\mathcal Inv^\sim(A)= \mathcal I(A)$. 

\begin{conjecture}\label{ConjU} The class $ \mathfrak I_*$ coincides with the class of all abelian groups, i.e., a subgroup $H$ of  an abelian group $A$ is uniformly fully inert if and only if $H$ is commensurable with a fully invariant subgroup of $A$. 
\end{conjecture} 

Let us note that the positive answer to the above conjecture implies the equality $\mathfrak I = \mathfrak I_u$, or
equivalently, the inclusion $\mathfrak I_u \subseteq \mathfrak I_*$.

We are not aware even if the following  weaker conjecture holds true: 

\begin{question} Is the equality $\mathfrak I = \mathfrak I_u$ $($or, equivalently, the  inclusion $\mathfrak I_u\subseteq \mathfrak I_*)$ true?
\end{question} 

\medskip
We shall see in the sequel that the known results support our stronger Conjecture \ref{ConjU}.
\medskip
 
Now we state a useful proposition which shows that the study of (uniformly) fully inert subgroups of finite direct sums can be reduced to box-like subgroups. It was proved in \cite[Lemma 4.1]{DGSV} for fully inert subgroups of {\em divisible} abelian groups, but the proof remains valid for (uniformly) fully inert subgroups of {\em arbitrary} abelian groups.

\begin{proposition}\label{Lemma:FI1}\cite{DGSV}  Let  $ A = A_1 \oplus\cdots\oplus A_n$ be an abelian group and $H\le A$. For $i=1,\ldots,n$, let $H_i:=H \cap A_i $ and   $H_* := H_1\oplus \cdots\oplus H_n$. Then the following conditions are equivalent:
\begin{itemize}
  \item[(a)] $H \in \mathcal I(A)$ $($resp., $H \in \mathcal I_u(A))$;
  \item[(b)] $[H:H_*]< \infty$  and $H_*\in \mathcal I(A)$ $($resp., $H_* \in \mathcal I_u(A))$. 
\end{itemize}
In case these conditions are fulfilled, then $H_i\in \mathcal I (A_i)$ $($resp., $H_i\in \mathcal I_u(A_i))$ for all $i$.
\end{proposition}

The next theorem from \cite{DGSV} gives some necessary conditions satisfied by a fully inert subgroup of a divisible abelian group. 

\begin{thm}\label{nec:condition}
\cite[Theorem 3.8]{DGSV} A proper fully inert subgroup $H$ of a divisible abelian group $D$ has finite torsion-free rank and $r_0(H) = r_0(D)$. 
Moreover, if $H$ has some non-trivial divisible $p$-quotient, then $H$ contains the p-primary component of $D$. 
\end{thm}

The equality $r_0(H) = r_0(D)$ is not given in \cite{DGSV}, but it easily follows, as the existence of a  non-torsion element $x\in D$ such that $\langle x \rangle \cap H = \{0\}$, will produce an endomorphism $\p\in \End(D)$ such that $x\in \p(H)$, so $r_0(H + \p(H))/H)>0$, witnessing $H \not \in   \mathcal I_\p(G)$, a contradiction. 

The above theorem and the next proposition provide easy examples showing that neither $\mathfrak I$, nor $\mathfrak I_u$ coincide with the class of all abelian groups.  

\begin{proposition}\label{new:lemma-}\hfill \begin{enumerate}
  \item Every finitely generated full rank subgroup of a torsion-free abelian group $A$ of finite rank is fully inert. 
  \item $ \{\{0\}, \Q^n\}= Inv (\Q^n) = Inv^\sim (\Q^n)= \mathcal I_u(\Q^n)$. 
  \item $H\in  \I(\Q^n)$ if and only if $H=H_1\oplus\ldots\oplus H_n$ where the $H_i$ are all isomorphic to a fixed subgroup of $\Q$.  
\end {enumerate}
Hence $\Q^n \not \in \mathfrak I_u$.
\end{proposition}

\pf (1) has been proved in  \cite[Lemma 2.4]{DGSV2} (see also \cite[Lemma 2.3]{DGSV}). 

To prove (2) let $\{0\}\ne H\in  \mathcal I_u(\Q^n)$. If $H$ is divisible, we have $r_0(H)=n$ and consequently $H=\Q^n$, by Theorem \ref{nec:condition}. If $H$ is not divisible, then there is a prime $p$ such that $H\ne pH$. The multiplication $\p={1\over p}$ is an automorphism of $A$ with $H<\p(H)$. Therefore, $|\p(H)/H|\ge p$ and -by induction- $|\p^n(H)/H|\ge p^n$ contradicting $H\in  \mathcal I_u(A)$.
This proves the equality   
$ \I_u (\Q^n)=\{\{0\}, \Q^n\}$ and, in view of (*), also (2). 

To prove (3), pick $H\in \mathcal I (\Q^n)$  
and note that $r_0(H) =n$, by Theorem \ref{nec:condition}. Now the assertion follows from  Proposition \ref{Lemma:FI1} and  \cite[Theorem 4.9]{DGSV}. \QED

%%%%%%%%%%%%%%%%

\subsection{Partial solutions of Problem \ref{Problem:I}}

Dikranjan, Salce and Zanardo  \cite{DSZ} have shown that the class $\mathfrak I$ contains all free abelian groups: 

\begin{thm} \cite{DSZ} A fully inert subgroup of a free abelian group $A$ is commensurable with a fully invariant subgroup of $A$. \end{thm}

\hskip-5mm Recall that a fully invariant subgroups  of a free abelian group $A$ must be of the form $nA$ for some $n\in\N$. 

In the case of periodic abelian groups, the problem is more complicated reflecting the fact that fully invariant subgroups of $p$-groups have a more complicated structure. On one hand, Goldsmith, Salce and Zanardo \cite{GSZ} have shown that $p$-groups with some further restriction  belong to $\mathfrak I$: 
 
\begin{thm} \cite[Theorem 3.10]{GSZ}  A fully inert subgroup of a direct sum of cyclic p-groups $A$ is commensurable with a fully invariant subgroup of $A$.
\end{thm}

On the other hand, one has the following example of a $p$-group that  does not belong to $\mathfrak I$. It is built using a construction of Pierce of abelian $p$-groups 
$A$ whose endomorphism ring  $End(A)$ is not that large, as $\End(A) = J_pid_A + E_s(A)$, where $J_pid_A=M(A)$ 
is the ring of multiplication by $p$-adics and $E_s(A)$ is the ideal of {\em small} endomorphisms $\p$ (i.e., such that for each  $k\in \N$ there is an $m\in \N$ with $\p((p^mA)[p^k]) = 0$). 

\begin{example}\cite{GSZ} 
There exist $p$-groups $A$ that contain fully inert subgroups which are not commensurable with any fully invariant subgroup of $A$, i.e., $A$ does not belong to $\mathfrak I$.
\end{example}

Recently also a result of Chekhlov  \cite{Ch}  about completely decomposable finite rank torsion-free abelian group has appeared.  Here we state an extension of that result.

\begin{thm} \cite{Ch} For a completely decomposable finite rank torsion-free abelian group $A = A_1\oplus \ldots \oplus A_n$, where  $r_0(A_i) = 1$ for each $i = 1,\ldots, n$, the following are equivalent:   
\begin{enumerate}
\item $A\in \mathfrak I$, i.e., every fully inert subgroup  is commensurable with a fully invariant subgroup;
\item  $pA_i \ne A_i$  and the types of $A_i$ and $A_j$ are either equal or incomparable, for each $i,j$ and each prime $p$.
\item $A\in\mathfrak I_u $ i.e., every fully inert subgroup is uniformly fully inert.
\end{enumerate}
 \end{thm}

\pf The equivalence (1) $\leftrightarrow$ (2) is proved in  \cite{Ch}. On the other hand, the implication (1) $\to$ (3) is trivial. 

\bigskip

To prove (3) $\to$ (2), assume by contradiction that $pA_1 = A_1$ for some prime $p$. Then the multiplication $\p_1$ by $1 \over p$ is an automorphism of $A_1$. Pick for every $i= 1,\ldots, n$ a cyclic subgroup $C_i$ of $A_i$ and let $H = C_1\oplus\ldots\oplus C_n$. This subgroup $H$ is fully inert, by Proposition \ref{new:lemma-}, while the subgroup $C_1 \cong \Z$ of $A_1$ is not uniformly fully inert since for each $m$ we have $|C_1+\p_1^{m}(C_1)/C_1|=p^m$. 
 Moreover, $\p_1$ may be extended to an endomorphism $\p$ of $A$ such that $
\p\restriction_{A_i}=0$ for each $i\ne 1$. Then $|H+\p^m(H)/H|\ge p^m$ for each $m$ and $H$ cannot be uniformly fully inert in $A$ either.

To prove the second part of the statement, argue by contradiction and assume (without loss of generality) that 
 the type of $A_1$ is strictly smaller that the type of $A_2$. 
Pick $x_i\in A_i$ such that for the $p$-heights we have $h_{p_m}(x_1)< h_{p_m}(x_2)$ for infinitely many primes $p_m$ with $m\in\N$. Then, for every $m\in \N$  there is a homomorphism $\psi_m: A_1 \to A_2$ such that  $\psi_m(x_1) = \frac{1}{p_1\cdot \ldots \cdot p_m}x_2$. Let $C_i:=\la x_i\ra$. Then  $C_2\leq \psi_m(C_1)$ and $|\psi_m(C_1):C_2| \geq m$. We can extend $\psi_m$ to an  endomorphism $\psi$ of the whole $A$ such that $\psi\restriction_{A_i}=0$ for each $i\ne 1$. Then we can show, as above, that  $\psi_m$ witnesses the failure of 
  $H = C_1\oplus\ldots\oplus C_n$ to be uniformly  fully inert, although $H$ is fully inert, a contradiction. \qed

%%%%%%%%%%%%%%%%%%%%%%
The characterization of the completely decomposable torsion-free abelian groups $A$  such that every fully inert subgroup of $A$ is commensurable with a fully invariant subgroup was extended in a more recent paper of Chekhlov  \cite[Theorem 2]{Ch+} to the case when $A$ has finitely many 
homogeneous components.  
%%%%%%%%%%%%%
\vskip8mm

\subsection{Self-inert groups and fully inert subgroups of divisible groups}

Although the case of divisible groups was the first one to be treated in the pioneering paper \cite{DGSV}, we are facing it only here as it is notoriously more complicated. We call an abelian group $G$  {\em self-inert}, if $G$ is fully inert in its divisible hull. This concept was introduced in \cite{DGSV}
under the term {\em inert}, that we avoid here for obvious reasons. 

\medskip

We collect in the next theorem the essential facts necessary for the description of the self-inert groups. 

\begin{thm}\label{looong:them}
\cite[Theorem 4.10]{DGSV} If $A$ is a self-inert group, then $A=t(A)\plus B$ splits and the subgroups $t(A)$ and $B$ are both self-inert.  
\begin{enumerate}
\item \cite[Theorem 4.9]{DGSV} A torsion-free abelian group $A$ is self-inert if and only if it is either divisible or completely decomposable homogeneous of finite rank.
\item \cite[Theorem 5.2]{DGSV} An abelian $p$-group $A$ is self-inert if and only if it is either divisible or bounded with at most one infinite Ulm-Kaplansky invariant.
\item \cite[Theorem 5.3]{DGSV} A torsion abelian group $A$ is self-inert if and only if each p-primary component of $A$  is a self-inert $p$-group and all but finitely many of them are either divisible or have a single nonzero Ulm-Kaplansky invariant.
\end{enumerate}
\end{thm}  

The conditions provided by the first assertion of  above theorem are only {\em necessary}. Conditions that are simultaneously 
{\em necessary and sufficient} for a mixed abelian group to be self-inert can be found in \cite[Theorem 5.5]{DGSV}. We skip them as they are rather technical and will not be used in the sequel in this survey. 
\medskip

The following simple observation can be helpful in the case for fully invariant subgroups: 

\begin{remark}\label{remark:sec5} {Let $\p$ an endomorphism of an abelian group $A$.}
\begin{enumerate}
  \item  If $H$ is a subgroup of an abelian group $A$ containing a fully invariant subgroup $N$ of $A$ and 
${\bar\p}: A/N \to A/N$ is the naturally induced endomorphism of $A/N$, then $H\in    \mathcal I_\p(A)$ if and only if $HN \in    \mathcal I_{\bar\p}(A/N)$. 
  \item Combining item (1) with Proposition \ref{new:lemma-} , we deduce that if $r_0(A) < \infty$ and $H$ is a subgroup of $A$ containing $t(A)$ with $r_0(H) = r_0(A)$, then  $H \in \mathcal I(A)$ if $H/t(A)$ is finitely generated. 
\end{enumerate}
 \end{remark}

Now we resolve Problem \ref{Problem:I} for divisible groups: 

\begin{thm}  For a divisible  abelian group $D$ the following are equivalent: 
\begin{itemize}
  \item[(1)] $D\in \mathfrak I$; i.e., every fully inert subgroup of $D$ is commensurable with a fully invariant subgroup of $D$; 
  \item[(2)] either  $r_0(D) = 0$ or  $r_0(D)$ is infinite;
  \item[(3)] $D\in \mathfrak I_u$, i.e., every fully inert subgroup of $D$ is uniformly inert. 
\end{itemize}
\end{thm}

\pf 
The implication $(1) \to (3)$ is trivial. In order to prove the implication $(3) \to (2)$ in a counter-positive form we need to check first that a non-torsion divisible group $D$ of finite torsion-free rank does not belong to $ \mathcal I_u$. Indeed, if $0<r(D)=n< \infty$, then we can write, up to isomorphism, $D = \Q^n \oplus t(D)$. Now the subgroup $G = \Z^n \oplus t(D) $ of $D$ is  fully inert by Remark \ref{remark:sec5}. On the other hand, 
$H$ is not uniformly fully inert in view of the last assertion of Proposition \ref{Lemma:FI1} and Lemma \ref{new:lemma-}.

To prove  the remaining implication $(2) \to (1)$ we note first that  according to Theorem \ref{nec:condition}, a divisible group $D$ of infinite torsion-free rank admits no proper fully inert subgroups, hence $D \in \mathcal I$ vacuously. Therefore, we are only left with the proof of the fact that all torsion divisible abelian groups $D$ belong to $\mathcal I$. 

Assume that $G\in \mathcal I(D)$, we have to prove that $G \sim H$ for some fully   invariant subgroup $H$ of $D$. Let $D = \bigoplus_p D_p$ and $G = \bigoplus_p G_p$ be the respective primary decompositions of $D$ and $G$.  Clearly, $G\in \mathcal I(D)$ if and only if $G_p\in \mathcal I(D_p)$ for every $p$ and $G_p$ is a fully invariant subgroup of $D_p$ for all but finitely many $p$.  Let 
$$
\Pi =\{p\in \mathbb P: \mbox{$G_p$ is a fully invariant subgroup of $D_p$}\}\  \   \mbox{ and }\  \  \pi = \mathbb P \setminus \Pi.
$$ 
Then $G_1 =  \bigoplus_{p\in \Pi} G_p$ is a fully invariant subgroup of $D$ and $G = G_1\oplus G_\pi$, where $G_\pi =  \bigoplus_{p\in \pi} G_p$.  Our next aim (see Claim) will be to prove that for every $p\in \pi$ there exists a  fully invariant subgroup $H_p$ of $D_p$ such that $G_p \sim H_p$ for every  $p\in \pi$. Then $H_\pi :=  \bigoplus _{p\in \pi} H_p$ is a fully invariant subgroup of $ \bigoplus _{p\in \pi} D_p$ and  $G_\pi \sim H_\pi$, so $G \sim G_1 \oplus H_\pi$, and the latter is a   fully invariant subgroup of $D$. 

\medskip

\noindent {\bf Claim. } For every $p\in \pi$  there exists a  fully invariant subgroup $H_p$ of $D_p$ such that $G_p \sim H_p$.

\pf If $G_p = D_p$, there is noting to prove. Hence, assume that $G_p \ne D_p$ in the sequel. 

It will be convenient to fix the divisible hull $\bar G_p$ of $G_p$ in $D_p$ as a minimal divisible subgroup of $D_p$ containing $G_p$. Then $D_p= \bar G_p\oplus D^*$ splits.  This yields that $G_p\in  \mathcal I(\bar G_p)$, i.e., $G_p$ is self-inert. By item (2) of Theorem \ref{looong:them}, $G_p$ is bounded and has at most one infinite Ulm-Kaplansky invariant.  Since $G_p$ is an essential subgroup of $\bar G_p$, this means that $G_p \sim H:= \bar G_p[p^m]$ for some $m \in \N$. If $r_p(G_p) < \infty$, then $G_p$ and $H$ are finite, so $G_p \sim \{0\}$, and we are done. This proves the claim. 
 
 \medskip

 Assume that $r_p(G_p)$ is infinite. Then $r_p(D^*) < \infty$, since otherwise one can easily   build an endomorphism of $D_p$ witnessing the failure of $G_p$ to be fully inert in $G_p$ (by simply sending the $p$-socle $G_p[p]$ to an infinite subgroup of the $p$-socle  of $D^*[p]$ and then extending this to an endomorphism of the divisible group $D_p$).  By $r_p(D^*) < \infty$, we deduce that $D^*[p^m]$ is finite, so $H\sim H \oplus D^*[p^m] = D[p^m]$.   This proves the claim (and the theorem), as $G_p \sim H$.  \QED

The next corollary confirms Conjecture \ref{ConjU} in the case of divisible groups: 

\begin{corollary} 
The class $\mathfrak I_*$ contains all divisible abelian groups, i.e., every uniformly fully inert subgroup of a divisible abelian group $D$ is commensurable with a fully invariant subgroup of $D$. 
In particular, $ Inv^\sim(t(D)) \cup \{D\} = Inv^\sim(D)$ and 
\begin{enumerate}
   \item if \ $0<r_0(D)<\infty$, then  $\mathcal I(D) \ne \mathcal I_u(D) = Inv^\sim(D)$, hence $D \in \mathfrak J_*$; 
   \item otherwise,  $\mathcal I(D) = \mathcal I_u(D)= Inv^\sim(D)$, hence $D \in \mathfrak I$.
\end{enumerate}
\end{corollary} 

\pf Let $D$ be a divisible group. The equality $ Inv^\sim(t(D)) \cup \{D\} = Inv^\sim(D)$ is obvious. 
 If $D\in \mathfrak I$, then $\mathcal I(D) = Inv^\sim(D)$ and there is nothing to prove. Assume that $D\not\in \mathfrak I$. Then $D = \Q^n \oplus t(D)$ for some positive $n \in \N$, by the above theorem. Let $G\in \mathcal I_u(D)$. We have to prove that $G$ is commensurable with a fully invariant subgroup. If $G\leq t(G)$, this follows by the above theorem. Let us see that this is the only case that may occur, in case $G\ne D$. 
Assume that $G \ne t(G)$. By Lemma \ref{Lemma:FI1}, $G \sim G_1 \oplus G_2$, where $G_1 \in \mathcal I_u(\Q^n)$ and  $G_2 \in \mathcal I_u(t(D))$. Our assumption gives $G_1\ne 0$. By Proposition \ref{new:lemma-},  from $G_1 \in \mathcal I_u(\Q^n)$ we can deduce that $G_1 = \Q^n$. From Theorem \ref{nec:condition} we deduce that $G_2 = t(D)$, i.e., $G = D$. 

%%%%%%%%%%%%%%%

\subsection{Generalizations}

To overcame the problem that an abelian group  may have  many fully invariant subgroups, recently,  Calugareanu \cite{Cal} has introduced    
{\it strongly invariant} subgroups, that is  subgroups $H$ of a (possibly non-abelian) group $G$ with the property that $f(H)\le H$ for each homomorphism $f:H\to G$. For example, fully invariant direct factors are clearly strongly invariant. On the other hand if $\Q_p = \{ m/n\in \Q\ |\ gcd(n; p) = 1\}$, then  the subgroup $p\Q_p$ of $\Q_p$ is fully invariant but not a strongly  invariant subgroup of $\Q_p$ (consider multiplication by $1/p$).

 Notice that unfortunately the usage of word  ``strong" in \cite{Cal}  is different from that in \cite{DGMT}.  
\begin{thm}\label{viva:Roumania}\cite{Cal}
The only strongly invariant subgroups of a reduced abelian $p$-group $A$ are the subgroups $A[p^n]$ for natural numbers n.
\end{thm}

Following Braez and Calugareanu \cite{BC}, in   the next statement we call a subgroup $H$ of an  abelian group $A$ {\it strongly inert}, if the factor group $ (H +f(H))/H$ is finite for every homomorphism $f: H\to A$. In the vain of Theorem \ref{viva:Roumania} we have:

\begin{thm}\cite{BC} Let $H$ a  strongly inert  subgroup (in the terminology of \cite{BC}) of an abelian group $A$. 
Then:
\begin{enumerate}
\item if $A$ is free abelian, then $A/H$ is finite.
\item if $A$ is torsion, then $H$ is commensurable with a strongly invariant subgroup of $A$.
\end{enumerate}
\end{thm}

Inertiality properties for modules have been considered, along the same lines as for abelian groups,  by Goldsmith, Salce and Zanardo \cite{GSZ-J} in the case of a $J_{p}$-module $A$, where $J_{p}$ denotes the ring of $p$-adic integers. A submodule $H$ of $A$ is called fully inert if $(H+\varphi(H))/H$ is finite for all $\varphi\in {\rm End}_{J_{p}}(A)$. Two submodules $K,L $ of $A$ are called commensurable if they are commensurable as subgroups, i.e., if $(K+L)/(K\cap L)$ is finite. We denote by $\mathfrak I_p$ the class of all $J_p$-modules $A$ such that every fully inert submodule of $N$ is commensurable with a fully invariant submodule of $A$. 

\begin{example}\cite[Proposition 2.3]{GSZ-J} If $A$ is a free $J_{p}$-module of finite rank $n$, then one can easily see that a non-zero fully invariant submodule of $A$
has the form $p^mA$ for some non-negative $m\in \Z$ and $p^mA\sim A$. Furthermore, a non-zero submodule $H$ of $A$ is fully inert if and only if 
$H$ has rank $n$. Clearly, this occurs precisely when $H$ has finite index (equivalently, when $A/H$ is torsion). In such a case, $H$ is commensurable with $A$. Therefore, $A \in \mathfrak I_p$.
\end{example}

These authors showed that the essential part of this property remains true for all free $J_{p}$-modules: 

\begin{thm}\label{J_p}\cite{GSZ-J}
 Let $A$ be a free $J_{p}$-module of infinite rank. Then for a submodule $H$ of $A$ the following are equivalent: 
  \begin{itemize} 
  \item[(1)] $H$  is fully inert;  
  \item[(2)] $H$ is commensurable with some fully invariant submodule; 
  \item[(3)] $A/H$ is a bounded $p$-group with at most one infinite Ulm-Kaplansky invariant and $(H+pA)/pA$ is finite whenever $A/H$ is infinite. 
\end{itemize}
\end{thm}

This theorem unifies Theorems 2.4 and 2.5 from \cite{GSZ-J}. The more technical condition in item (3) describes intrinsically the fully inert
 submodules. The assumption $|(H+pA)/pA|<\infty$ cannot be avoided in general, as shown in \cite[Example 2.6]{GSZ-J}.

The equivalence of (1) and (2) in Theorem \ref{J_p} and the next theorem show that free $J_{p}$-module of infinite rank as well as complete torsion-free  $J_{p}$-modules belong to the class $\mathfrak I_p$.  

\begin{thm}  \cite[Theorem 3.6]{GSZ-J} A submodule of a complete torsion-free  $J_{p}$-module $A$ 
is fully inert if and only if it is commensurable with a fully invariant submodule (that is with some $p^nA$ for $n\in \N$).
\end{thm}

Finally, an example of a torsion-free $J_{p}$-module $A$ with a fully inert submodule which is not commensurable with any of its submodules $p^{n}A$
is obtained using deep realization theorems of commutative rings as endomorphism rings of $J_p$-modules (\cite[Theorem 4.1]{GSZ-J})

%______INIZIO SEZIONE_6___________________

\section{Applications to algebraic entropy of group endomorphisms}\label{appl} 

%\pagebreak
We see in this section that the concept of $\p$-inert subgroup is a very helpful tool also for the study of the dynamical properties of a given endomorphism. 

\subsection{Algebraic entropy} \label{entropy}

Roughly speaking, the algebraic entropy is a non-negative real number or $\infty$ assigned to an endomorphism $\p$ of a group $G$ 
measuring how chaotically $\p$ acts on certain families of (most often, finite) subgroups or subsets of $G$.

Since the entropy is defined as a limit, one needs a tool providing the existence of limits.  

\begin{remark}\label{subadditive}
A sequence $\{a_n\}_{n\in\N}$ of non negative real numbers is called {\em subadditive} if  $a_{n+m}\leq a_n  +a_m$ for every $n,m\in \N$. 
According to a folklore fact, known as {\em Fekete Lemma}, for a subadditive sequence $\{a_n\}_{n\in\N}$, the sequence $\{\frac{a_n}{n}\}_{n\in\N}$ converges and $\lim_{n\to\infty}\frac{a_n}{n}=\inf_{n\in\N}\frac{a_n}{n}.$
\end{remark}

\subsubsection{The entropies $\ent$ and $h_{alg}$ and their basic properties}

For a subset $F$  of an abelian group $(A,+)$, $\p\in End(A)$  and  $n\in \N$,  the {\em $n$-th trajectory} of $F$ under the action of $\p$ is the set  
$$
T_n=T_n(\p,F):=F+\p(F)+\ldots+\p^{n-1}(F).
$$ 
If $F$ is a subgroup, then each $T_n:=T_n(\p,F)$ is a subgroup and the {\em trajectory} $T(\p,F):= \bigcup_n T_n$ of $F$ under the action of $\p$ is the smallest $\p$-invariant subgroup of $A$ containing $F$. The entropy measures how rapidly grows the chain of $n$-th trajectories $T_n(\p,F)$ when approximating the trajectory $T(\p,F)$. 

When $F$ is finite, the limit 
\begin{equation*}
\mathcal 
H_{alg}(\p,F):= \lim_{n\to\infty}\frac{\log|T_n(\p,F)|}{n}
\end{equation*}
exists, according to Remark \ref{subadditive}, as the sequence $\log|T_n(\p,F)|$ is subadditive. This idea was originally used by  Adler, Konheim and McAndrew \cite{AKM} in the case when $F$ is a finite {\em subgroup}, so they defined the \emph{algebraic entropy} $\ent$ of $\p$ by  
\begin{equation*}
\ent(\p):=\sup\{\mathcal H_{alg}(\p,F): F\ \text{finite subgroup of}\ A\}.
\end{equation*} 

\begin{example}[Bernoulli shift]\label{Ber}
The leading example is the (right) Bernoulli shift defined as follows. Let $F$ be a finite abelian group and let 
$\beta_F : \bigoplus_{n=1}^\infty F \to \bigoplus _{n=1}^\infty F$ be the endomorphsm defined by 
$$
\beta_F(x_1, x_2, x_3, \ldots ) = (0, x_1, x_2, x_3, \ldots ).
$$
Then $\ent(\beta_F) = \log |F|$ \cite{DGSZ}. The Bernoulli shift $\beta_F$ can be defined for an infinite $F$ as well, but then $\ent(\beta_F) = \infty$. 
\end{example}

The entropy function $\ent$ is not appropriate in studying endomorphisms of torsion-free groups. To this end one can define an appropriate modification as follows. The \emph{algebraic entropy} $h_{alg}(\p)$ of $\p$ is defined as   
\begin{equation}\label{def:h}
h_{alg}(\p):=\sup\{\mathcal H_{alg}(\p,F): F\ \text{finite subset of}\ A\}.
\end{equation}

 The use of the same term for two distinct entropy functions may look somewhat embarrassing, but we are following the traditional 
terminology. Anyway, if will always be clear from the context which is the precise entropy in question.  

The subtle question of whether the $\sup$ in (\ref{def:h}) is attained (i.e., $\sup$ is a $\max$) in case when $h_{alg}(\p)<\infty$ is still open and will
be discussed in Remark \ref{realiz}.  Clearly, $\ent(\p) = h_{alg}(\p\restriction_{t(A)})$ for any abelian group $A$ and $\p\in \End(A)$, so  $\ent(\p) \leq h_{alg}(\p)$
and these entropies coincide when $A$ is torsion. On the other hand, if $A$ is not torsion, then there exists $\p\in \End(A)$ with $\ent(\p) < h_{alg}(\p)$
(any multiplication by an integer $m>1$ will do, see Example \ref{AYF}). 

One can prove that $h_{alg}(\p^{-1}) = h_{alg}(\p)$ for an automorphism $\p$, and more generally 
$$
h_{alg}(\p^{n}) = |n| h_{alg}(\p) \mbox{ for every integer }n.\eqno(LL)
$$ 
This rule for computation of the entropy of powers is known as Logarithmic Law. It remains true also for endomorphisms as far as $n > 0$.  
Note that (LL) implies $h_{alg}(id) = 0$, although this can be directly obtained from the definition. We shall see below that this need not be true in the non-abelian case. 

The above definition was proposed in \cite{DG2}, inspired by a similar notion of entropy introduced by Peters \cite{Pet}  in 1979 for automorphisms $\p$
and using the ``negative" $n$-th trajectories (i.e., $T_n(\p^{-1},F)$). Since the limit  obviously gives $h_{alg}(\p^{-1})$, 
%these 
the values of these two notions coincide for  automorphisms, by virtue of (LL). Another leading example of computation of the algebraic entropy is given below (Example \ref{AYF}).

In order to describe the algebraic entropy of endomorphisms of $\Q^n$ we need to recall the definition of Mahler measure. 
Let $n$ be a positive integer, let $f(t)=st^n+a_{n-1}t^{n-1}+\ldots+a_0\in\mathbb Z[t]$ be a primitive polynomial and let $\{\lambda_i:i=1,\ldots,n\}\subseteq\mathbb C$ be the set of all roots of $f(t)$ taken with their multiplicity. 
The (logarithmic) Mahler measure of $f(t)$ was defined independently by Lehmer  
and Mahler 
 in two different equivalent forms as  
$${m}(f(t))=\log|s|+\sum_{|\lambda_i|>1}\log|\lambda_i|.$$

\begin{example}[Algebraic Youzvinski formula]\label{AYF}
 An endomorphism $\p$ of $\Q^n$ is uniquely determined by its $n\times n$ matrix over $\Q$ and its characteristic polynomial $f(t) = st^n+a_{n-1}t^{n-1}+\ldots+a_0$ considered as a primitive polynomial with integer coefficients. The  algebraic entropy of $\p$ con be computed 
by using the formula $h_{alg}(\p) = m(f)$, established in \cite{GV} under the name {\em Algebraic Youzvinski formula}.   

   In particular, if  $\p$ is a non-zero endomorphism of $\Q$, then $\p$ is the multiplication by some rational number $a/b$, with non-zero co-prime $a,b$, so the above formula gives $h_{alg}(\p) =\max\{\log |a|, \log |b|\}$. 
 Therefore,  for every non-zero endomorphism $\p: \Z \to \Z$ one has $h_{alg}(\p) = \log |\p(1)|$. 
\end{example}

The next theorem collects three fundamental properties of the algebraic entropy:  

\begin{thm}\label{UT}\cite{DG2} Let $\p$ be an endomorphism of an abelian group  $A$.
\begin{itemize} 
\item[(1)]  {\em [Invariance under conjugation]} If $A_1$ is a group and $\xi : A \to A_1$ an  isomorphism, then for the endomorphism $\psi =\xi\p\xi^{-1} $ of $A_1$ one has $h_{alg}(\p)=h_{alg}(\psi)$.
\item[(2)] {\em [Continuity]} If $A$ is direct limit of $\p$-invariant subgroups $\{A_i:i\in I\}$, then $h_{alg}(\p)=\sup_{i\in I}h_{alg}(\p\restriction_{A_i}).$
\item[(3)] {\em [Addition Theorem]} If $N$ a $\p$-invariant subgroup of $A$ and $\overline\p:A/N\to A/N$ the endomorphism induced by $\p$. Then $$h_{alg}(\p)=h_{alg}(\p\restriction_N)+ h_{alg}(\overline\p).$$
\end{itemize}
\end{thm}

\subsubsection{Algebraic $i$-entropy for module endomorphisms}

Let $R$ be a ring and let $\Mod_R$ be the category of all left $R$-modules. An {\em invariant} of $\Mod_R$ is a function  $i: \Mod_R \to \R_{\geq0}\cup \{ \infty\}$ with $i(0)=0$ and $i(M) = i(N)$ whenever $M\cong N$. The invariant is said to be  \emph{subadditive}, if 
\begin{itemize}
  \item[(a)] $i(N_1 + N_2)\leq i(N_1) + i(N_2)$ for all submodules $N_1$, $N_2$ of $M$;
  \item[(b)] $i(M/N)\leq i(M)$ for every submodule $N$ of $M$;
\end{itemize}
Furthermore, the invariant $i$ is called {\em discrete}, if the set $I$ of all finite values $i(M)$ is a discrete subset of $\R$ (i.e., $I\cap [0,N]$ is finite for every $N >0$), integer valued  invariant are typical examples of discrete  invariants. 
The invariant $i$ is called \emph{additive} if equality holds in (a) (this implies (b)). In case $i$ is additive and 
$$
i(M) = \sup\{i(N): N\mbox{ is a finitely generated submodule of  }M\}
$$
the invariant $i$ is call a {\em length function}. 

\smallskip
For a right $R$-module $M$, denote by  $\sF_i(M)$ the family of all submodules $N$ of $M$ with $ i(N)< \infty$. The {\em algebraic $i$-entropy} $\mbox{ent}_i$ was introduced by Salce and Zanardo \cite{SZ} as follows.  

\begin{definition}{\rm \cite{SZ}}
Let $i$ be a subadditive invariant, $M$ be an $R$-module and $\p:M\to M$ an endomorphism.  The \emph{algebraic $i$-entropy of $\p$ with respect to $N\in\mathcal F_i(M)$} is 
\begin{equation*}
 \mathcal H_i(\p, N):=\lim_{n\to \infty}\frac{i(T_n(\p,N))}{n};
\end{equation*}
The \emph{algebraic $i$-entropy} of $\p$ is $\mbox{ent}_i(\p):=\sup\{H_i(\p,N): N\in \mathcal F_i(M)\}$.
\end{definition}

Typical examples are: 
\begin{itemize}
\item[(a)] $R=\Z$, so $\Mod_R$ is the category of abelian groups, now we have two examples: 
\begin{itemize}
  \item[(a$_1$)]  $i(M) = \log |M|$, then  $\ent_i$ coincides with $\ent$;  
  \item[(a$_2$)] $i(M) = r_0(M)$, this entropy, termed {\em rank-entropy}, was largely studied in \cite{SZ,GS};  
\end{itemize}
\item[(b)] $R$ is a field, so that $\Mod_R$ is simply the category of $R$-linear spaces, and $i(M) = \dim M$. This entropy $\ent_{\dim}$ was briefly mentioned in \cite{SZ} and thoroughly studied in \cite{Car,GBS}.
\end{itemize}

Salce, V\' amos and Virili \cite{SZV} extensively studied the entropy function $\ent_i$ for a lenght function $i$. 

\subsection{Intrinsic algebraic entropy} %\label{entropy}
Regardless of their nice properties, the entropies $\ent$ and  $h_{alg}$ have their disadvantages. The case of $\ent$ was discussed above,  the case of $h_{alg}$ has the obvious inconvenience to work with {\em subsets} rather than {\em subgroups}. An excellent remedy to both problems is provided by the
{\em intrinsic algebraic entropy} defined below.  
    
It was observed already in \cite{DGSZ}, that if $H$ is a subgroup of an abelian group $A$, then for the increasing chain 
of subgroups $T_n:=T_n(\p,H)$ of $A$ each quotient $T_{n+1}/T_n$ is an epic image of $H+\p(H)/H = T_2/T_1$ for every $n\in \N$. Therefore, when $H$ is finite, the sequence $|T_{n+1}/T_n|$  stabilizes to a value  that coincides with $\mathcal H_{alg}(\p, H)$. Clearly this occurs also when simply the quotient $H+\p(H)/H$ is finite, i.e., $H$ is $\p$-inert. This simple but relevant observation allowed Dikranian, Giordano Bruno, Salce and Virili \cite{DGSV2} to extend the definition of the entropy  $\mathcal H_{alg}(\p, F)$  \ to $\p$-inert subgroups $F$ of $A$.  The {\em intrinsic algebraic entropy}  is defined as
\begin{equation}\label{intri}
\widetilde{\ent}(\p) := \sup\{\mathcal H_{alg}(\p, F): F \mbox{ is a $\p$-inert subgroup of } A\}. 
\end{equation}
This dynamical invariant is better behaved since it allows to treat also endomorphisms of torsion-free abelian groups where the entropy $\ent$ vanishes
completely for the lack of non-trivial finite subgroups. One can easily show that $$ {\ent}(\p)\le \widetilde {\ent}(\p)\leq h_{alg}(\p),$$ where the first inequality becomes an equality when $A$ is torsion (\cite{DGSV2}).  The function  $ \widetilde {\ent}$ retains many of the useful properties of $\ent$: the Logarithmic Law (LL) for $n>0$ \cite[Lemma 3.12]{DGSV2}, all three properties (1)--(3) from Theorem \ref{UT} (these are Lemma 3.9, Lemma 3.14 and Corollary 5.10 respectively, from \cite{DGSV2}), etc.

Here we recall the following counterpart of the Algebraic Youzvinski formula obtained in \cite{DGSV2}. 

\begin{thm}\emph{\cite[Theorem 4.2]{DGSV2}}\label{TY-intr}
Let $m\in \N$ and let $\p:\Q^m\to \Q^m$ be an endomorphism. Then $\widetilde{\ent}(\p) = \log s$, where $s$ is the positive leading coefficient of the characteristic primitive polynomial of $\p$ over $\Z$. 
\end{thm}

Further details and results can be found in \cite{DGSV2} and in the survey paper \cite{GSs}.

The part {\em intrinsic} in the name of this new entropy is justified by the fact that in the defining formula 
(\ref{intri}) one makes recourse to a family of subgroups of $A$ depending (intrinsically) on $\p$. 
Recently, Goldsmith and Salce \cite{GS} investigated the ground to consider this entropy really as an {\em intrinsic} one. To this end they 
introduced the so called {\em fully inert algebraic entropy} 
$$
 \tilde{\ent}(\p) = \sup\{\mathcal H(\p, F): F\in \mathcal I(A)\},
  $$ 
  where $A$ is an abelian group and $\p\in \End(A)$. In the definition of this entropy the choice of the subgroups $F$ does not depend on $\p$. 
 Their aim was to compare this new entropy with the intrinsic one.     They showed that $ \widetilde{\ent}(\p)$ can be computed by using the fully inert algebraic entropy $
 \tilde{\ent}(\bar \p)$ of the induced endomorphism $\bar \p$ of  $\p$-invariant sections of $G$ \cite[Theorem 5.5]{GS}. On the other hand, it is also shown that  the fully inert algebraic entropy need not coincide with the intrinsic one in general and cases when they coincide (so ``when the intrinsic entropy is not really intrisic") are examined (\cite[Theorems 5.1, 5.3 and 6.4]{GS}).

\subsubsection{Intrinsic algebraic $i$-entropy}\label{i-i-entropy}

Assume that $i$ is a discrete invariant of $R$-modules. Then one can define  intrinsic algebraic $i$-entropy for endomorphisms of 
$R$-modules following the pattern of the definition of   intrinsic algebraic entropy for abelian groups. Namely, noticing that if $F$ is a submodule of $M$, then each $T_n:=T_n(\p,F)$ is a submodule of $M$ and each quotient $T_{n+1}/T_n$ is an epic image of $T_{n}/T_{n-1}$, so in particular of 
$F+\p(F)/F = T_2/T_1$. Therefore, when $i(F+\p(F)/F)$ is finite, the sequence $i(T_{n+1}/T_n)$  stabilizes to a value, say  $\mathcal H_i(\p, F)$. This motivates the notion of an {\em $i$-$\p$-inert submodule } of $M$, namely a submodule $F$ of $M$ such that $i(F+\p(F)/F)$ is finite. Then the {\em intrinsic algebraic $i$-entropy} $\widetilde{\ent}_i$ is defined by 
$$
\widetilde{\ent}_i(\p) = \sup\{{\mathcal H}_i(\p, F): F \mbox{ is a $i$-$\p$-inert submodule of } M\}.
$$

In the sequel we shall use it for linear spaces and $i = \dim$. In this setting the intrinsic algebraic $i$-entropy was defined and studied in \cite{Car}
(see also the recent paper \cite{CGB} and the comment in \S \ref{Topo} below). It was proved in \cite{Car}, that for linear space endomorphisms the entropies $\ent_{\dim}$ and 
$\widetilde{\ent}_{\dim}$ coincide.
The linear maps $\p$ for which each subspace is $\dim$-$\p$-inert have been treated in \S\ref{linear}.

%%%%%%%%%%%%%%%%%%%%%%%%%%%%%%%%%%%%%%%%%%%%%%%%

\subsection{Values and uniqueness of the algebraic entropy} \label{values:entropy}

It follows from the definiton of $\ent$ and $\widetilde{\ent}$, that the finite values of these entropies belong to the (discrete) subset $\log \N = \{\log n: n \in \N, n>0\}$
of the whole possible range $\R_{\geq 0}$.

As far as the entropy function $h_{alg}$ is concerned, the question of understanding its precise range, i.e., set of values 
$$
\mathcal E =\{h_{alg}(\p): \p \in End(A), A \mbox{ abelian group}\}\subseteq \R_{\geq 0} \cup \{\infty\}
$$
is much more complicated and still open. In fact, it is not known whether $\mathcal E= \R_{\geq 0} \cup \{\infty\}$. More precisely, even the simpler question of whether  $0= \inf (\mathcal E\setminus \{0\})$ is open and equivalent to a long standing open problem of Lehmer (see Problem \ref{LP}). 

The Lehmer Problem concerns the range of the  Mahler measure. Let us start by pointing out that  the case of zero Mahler measure is completely determined by  the a
classical  theorem due to Kronecker: $m(f(t))=0$ for a primitive $f(t)\in\Z[t]\setminus\{0\}$ if and only if $f(t)$ is cyclotomic (i.e., all the roots of $f(t)$ are roots of unity).
Before going further let us recall, that Mahler measure is used also when referring to an algebraic integer  $\alpha$, taking into account the Mahler measure
of the minimal polynomial of $\alpha$. In these terms, the (still open) problem posed by Lehmer in 1933 aks whether the lower bound of all Mahler measures
of algebraic integers that are not roots of unity is positive:

\begin{problem}[Lehmer Problem]\label{LP} 
Given any $\delta>0$, is there any algebraic integer whose Mahler measure is strictly between $0$ and $\delta$? 
\end{problem}

One can show that the positive solution of Lehmer Problem is equivalent to the fact that $\mathcal E = \R_{\geq 0} \cup \{\infty\}$ (as well as to the apparently much 
weaker property $0\in \inf (\mathcal E\setminus \{0\})$, see \cite{DG*} and \cite[Theorem 1.5]{DG2}). 

The negative solution of Lehmer Problem has also a notable impact on the algebraic entropy (beyond the bold negation $\mathcal E\ne  \R_{\geq 0} \cup \{\infty\}$) 

\begin{remark}\label{realiz}
Let $A$ be an abelian group and $\p\in \End(A)$ with $h_{alg}(\p)<\infty$. 
We shall briefly say that $h_{alg}(\p)<\infty$ {\em realizes}, if there exists a finite subset $F$ of $G$ such that $h_{alg}(\p) = h_{alg}(\p,F)$.
It was shown in \cite{DG2} that the negative solution of Lehmer Problem is equivalent to the assertion ``the algebraic entropy of 
arbitrary endomorphisms of abelian groups realizes whenever it is finite".
\end{remark}

The importance of the  Algebraic Youzvinski formula (see Example \ref{AYF}) and the Bernoulli shifts becomes clear from the Uniqueness Theorem for the algebraic entropy \cite[Theorem 1.3]{DG2}, asserting that the algebraic entropy  $h_{alg}$ is the unique $\R_{\geq 0} \cup \{\infty\}$-valued function defined 
on endomorphisms of abelian groups having the three properties (1)--(3) from Theorem \ref{UT}, satisfying the Algebraic Yuzvinski formula and the ``normalization" condition, imposed by the values taken on the Bernoulli shifts $\beta_F$ (see Example \ref{Ber}). Similar uniqueness theorem for $\widetilde{\ent}$ is proved in \cite[Theorem 6.1]{DGSV2}. Both uniqueness theorems were inspired by the uniqueness theorem for the entropy function 
$\ent$ for endomorphisms of torsion abelian groups, previously established in \cite{DGSZ}. 

\subsection{Growth function of an endomorphism}\label{Intrinsic_non_abelian} 

Dikranjan and Giordano Bruno in \cite{DG*,DG:PC} have proposed an extension  of the notions of entropy in the non-abelian context  as follows.

If $G$ is any group in multiplicative notation, $F$ a finite subset of $G$ and $\p$ an  endomorphism of $G$, one can consider the trajectory $T_n:=T_n(\p,F):=F\cdot F^\p\cdot\ldots\cdot F^{\p^n}$ as the  setwise product of the images $F^{\p^i}$ for $i\le n\in \N$. Clearly $|T_n|\le |F|^n$. 
Furthermore, the limit ${\mathcal H}(\p,F):=\lim_n \frac{\log |T_n|}{n}$ exists by Remark \ref{subadditive} as the sequence $\log |T_n|$ is again subadditive 
 and the entropy $h_{alg}(\p)$ is the sup of ${\mathcal H}(\p,F)$ when $F$ ranges over all finite subsets of $G$. Now properties (1) and (2) from Theorem \ref{UT}
remain true, but (c) fails even for the identity automorphism, which may have infinite entropy in the non-abelian case. Thus, the Logarithmic Law (LL)
fails for $n=0$, but remains true for all $n\ne 0$ \cite{DG*}. Yet, it implies $h_{alg}(id_G) \in \{0, \infty\}$.

The above described ``anomalous" behavior of the identity automorphism in the non-abelian case can be understood under a more precise looking glass.
 Namely, considering the {\em growth function} $\gamma_{\p,F}(n)=|T_n(\p,F)|$ of the $n$-th trajectory rather than the limit ${\mathcal H}(\p,F):=\lim_n \frac{\log |T_n|}{n}$  which gives less information (see Lemma \ref{easy:lemma}). Following \cite{DG3,DG*}, we say that $\p$ {\em has polynomial $($respectively, exponential$)$  growth} if $\gamma_{\p,F}(n)$ has such growth for all 
(respectively, for one) finite subsets $F$ of $G$. Otherwise, we say that $\p$ has intermediate growth. The next easy to check 
property (see \cite{DG3,DG*}) gives the basic relation between entropy and growth of an endomorphism.

\begin{lemma}\label{easy:lemma}  Let $\p$ be an endomorphism of a group $G$
and let $F$ be a finite subgroup of $G$. Then  ${\mathcal H}(\p,F)>0$ if and only if $\gamma_{\p,F}(n)$ has exponential growth. In particular, 
 $h_{alg}(\p)>0$  if and only if $\p$ has exponential growth. 
\end{lemma}

  The above terminology meets the existing one on growth of groups in the following sense. A group $G$ has polynomial (respectively, exponential) growth (in the classical sense of geometric group theory \cite{dH}) precisely when the identity automorphism of $G$ has polynomial (respectively, exponential) growth. Therefore, in a group $G$ of exponential growth (e.g., a free non-abelian group) one has $h_{alg}(id_G)= \infty$. 

  This trend was triggered by the following result recalling relevant well known facts on the growth of finitely generated groups (according to the celebrated theorem of Milnor-Wolf every finitely generated soluble group has either polynomial or exponential growth \cite{dH}). 
  
  \begin{thm} \cite{DG3}\label{Pinsker}
All endomorphisms of an abelian group have either polynomial or exponential growth. 
\end{thm}

\vskip-1mm
One can see that  {\em all inner automorphisms of a group have the same growth type} \cite{GBSp}.  Moreover, these authors gave also a counterpart of Theorem \ref{Pinsker} for  locally finite groups.

\begin{thm}
\cite{GBSp} {An endomorphism $\p$ of a locally finite group has polynomial growth if and only if the algebraic entropy $h(\p)$ is  zero}.
\end{thm}
\vskip-1mm
According to Lemma \ref{easy:lemma}, this theorem implies that an endomorphism of  a locally finite group has either polynomial or  exponential growth. 
A simultaneous extension of this theorem and Theorem \ref{Pinsker} to locally (soluble-by-finite) groups 
was obtained recently by Giordano Bruno and Spiga: 

\begin{thm}\cite{GBS+}{ An endomorphism of  a  locally (soluble-by-finite) group has either polynomial or  exponential growth.}
\end{thm}

\vskip-1mm
 One can define  {\em intrinsic}  entropy in the non-abelian case making use, again, of inert subgroups.  
The idea is that if $H$ is a $\p$-inert subgroup of a group $G$, then the index $|H^\p:(H\cap H^\p)|$ is just the minimal size of a subset $X$ of $H$ such that $H^\p=(H\cap H^\p)X^\p$ (setwise product), hence $H^\p\subseteq H X^\p$. 
This will enable us to extend the definition of algebraic entropy (in the non-abelian setting) to an {\em intrinsic} one by using the trajectories of finite transverals like $X^\p$ above. 
We formulate this in the following immediate proposition: 
\vskip-1mm

\begin{proposition}\label{UD} Let $\p$ be an endomorphism of a group $G$.
 If $H$ is a $\p$-inert subgroup of $G$ with $t:=|H^\p:(H\cap H^\p)|$, then for every $n\in \N$, the setwise product  
$T_n:=H\cdot H ^\p\ldots H^{\p^n}$ is contained in a setwise product $H Y_n$,  where $Y_n\subseteq G$ is finite. 
If $t_n$ is the smallest possible size of a  finite $Y_n \subseteq G$ with $T_n \subseteq N Y_n$, then $t_n \le t^n$ and 
\begin{equation}\label{33}
\widetilde{\mathcal H}(\p,H) := \limsup \frac{\log t_n}{n}
\end{equation}
exists, with $\widetilde{\mathcal H}(\p,H) \leq t$. 
\end{proposition}
\pf Fix a subset $F$ with $|F|\leq t$ and $H^\p\subseteq HF^\p$. Then, by transforming this relation by $\p^n$, one has $H^{\p^{n+1}}\subseteq H^{\p^n}F^{\p^{n+1}}$. Let us show that the claim holds with $Y_n:=F^{\p} \ldots F^{\p^{n}}$ arguing by induction. The case $n = 1$ is simply our hypothesis on $F$. 
 Assume the statement holds for $n$.   Now 
 $$
 T_{n+1}=T_{n}H^{\p^{n+1}}\subseteq T_n H^{\p^n}F^{\p^{n+1}}=T_n F^{\p^{n+1}}\subseteq  H F^{\p} \ldots F^{\p^{n}} F^{\p^{n+1}}.
 $$ 
 The inequality $|Y_n|\le t^n$ is trivial. This proves $t_n \le t^n$. The rest is obvious. 
\QED

Call the limit in (\ref{33}) {\em intrinsic entropy of $\p$ with respect to} $N$. In case $G$ is abelian, the $\limsup$ in (\ref{33}) is a limit: 
 
\begin{problem}
Under which condition the $\limsup$ in (\ref{33}) defining the intrinsic entropy $\widetilde{\mathcal H}(N,\p)$ is a limit ?   
\end{problem}

\subsubsection{Limit-free computation of entropy}
The next rule for limit-free computation of entropy was proved in \cite{DG-lf} in order to repair a wrong statement in \cite{Y}:  

\begin{thm}\cite[Algebraic Formula]{DG-lf} Let $G$ be a  locally finite group, $\p:G\to G$ an endomorphism and $F$ a finite normal subgroup of $G$ such that $\ker\p\cap T(\p,F)$ is finite. Then 
$$
{\mathcal H}_{alg}(\p,F) = \log \left |\frac{T(\p,F)}{\p(T(\p,F))} \right|-\log|\ker\p\cap T(\p,F)|.
$$
\end{thm}

It was proved in \cite[Lemma 4.1]{DG-lf} that the restraint $|\ker\p\cap T(\p,F)|<\infty$ can be omitted in the abelian case; namely $\ker\p\cap T(\p,F)$ is finite whenever  $G$ is a torsion abelian group.
\smallskip

The next corollary was formulated by Yuzvinski's \cite{Y} who omitted the condition for injectivity of the endomorphism:  

\begin{corollary}\label{Coro0:May21}\cite[Corollary 1.1]{DG-lf} Let $G$ be a locally finite group, $\p:G\to G$ an injective endomorphism and $F$ a finite normal subgroup of $G$.  Then
$$
{\mathcal H}_{alg}(\p,F)  = \log \left |\frac{T(\p,F)}{\p(T(\p,F))} \right|.
$$
\end{corollary}

\begin{corollary}\label{Coro:ottobre}\cite[Corollary 2.3]{DG-lf}
Let $G$ be a locally finite group, $\p:G\to G$ an endomorphism and $F$ a finite normal subgroup of $G$ such that $G=T(\p,F)$. If $\ker\p$ is finite, then
$$
{\mathcal H}_{alg}(\p,F) = \log |\mathrm{coker}\, \p|-\log|\ker\p|.
$$
In particular, $|\mathrm{coker}\,\p|\geq|\ker\p|$.
\end{corollary}

 For limit-free computation of intrinsic entropy see \cite{GV1}. 

 \subsection{Inert subgroups of topological groups}\label{Topo}

There are some circumstances when inert subgroups naturally appear in the framework of topological groups
and their continuous endomorphisms. In the sequel we denote by $C\End(G)$ the family of all continuous endomorphisms of a topological group $G$. 

First we recall the definition of algebraic entropy of a continuous endomorphism $\p$ of a locally compact group $G$. Let $\mu$ be a Haar measure of $G$ and let $\mathcal U_G$ be the local base at  $1_G$ formed by all compact  \nbd s of $1_G$. Then $T_n(\p,U)$ is a compact set of $G$, so $\mu(T_n(\p,U))<\infty$. Call 
$$
\mathcal H_{alg}(\p,U) := \limsup_n \frac{\log \mu(T_n(\p,U))}{n} <\infty
$$
{\em algebraic entropy} of $\p$ with respect to $U$ and let $h_{alg}(\p)= \sup \{H_{alg}(\p,U): U \in  \mathcal U_G\}$, {\em algebraic entropy} of $\p$. 
For a discrete $G$ this coincides with the definition given in \S 6.3, so using the same notation does not lead to confusion.  

One can show that $h_{alg}(\p)$ does not depend on the choice of the Haar measure $\mu$. In the case of abelian  locally compact groups, this definition is due to Virili \cite{Vi}, the above definition is from \cite{DG*}. 

Let $G$ be a locally compact totally disconnected group. It is a well known fact, due to van Dantzig \cite{vD}, that the  filter base $\mathcal B_G$ consisting of all open compact subgroups of $G$ 
forms a base of the \nbd \ filter $\mathcal U_G$ at the identity of $G$. Here comes a remarkable connection to inertiality: 
the ring $C\End(G)$ is a subring of the ring $\mathcal I_{\mathcal B_G} \End(G)$ of all $\mathcal B_G$-inertial endomorphisms of $G$, in other words: 
  
\begin{lemma}\label{laaast:lemma}
Every $U\in \mathcal B_G$ is $\p$-inert for every $\p\in C\End(G)$. 
\end{lemma}

Indeed, $\p(U)$ is compact, as a continuous image of the compact subgroup $U$, whereas $U\cap \p(U)$ is an open subgroup of the  compact group $\p(U)$. This entails  finiteness of the index $|\p(U):(U\cap \p(U)|$, so $U$ is  is $\p$-inert. 

We do not know if one can invert Lemma \ref{laaast:lemma}, i.e., does the equality $\mathcal I_{\mathcal B_G} \End(G) = C \End(G)$ hold true?  

\begin{question}
If $G$ be a locally compact totally disconnected group, are then all $\mathcal B_G$-inertial endomorphisms of $G$ necessarily continuous? 
\end{question}

One can apply Proposition \ref{UD} and Lemma \ref{laaast:lemma} to get a a measure-free formula for the algebraic entropy in totally disconnected locally compact groups. 

\begin{proposition}\label{LAST}
Let $G$ be a totally disconnected locally compact group, $\p:G\to G$ a continuous endomorphism and $U\in\mathcal B(G)$. Then
${\mathcal H}_{alg}(\p,U) = \widetilde{\mathcal H}(\p,U)$.  
\end{proposition}

\pf 
By Lemma \ref{laaast:lemma}, $U$ is $\p$-inert so by Proposition \ref{UD} one can define $\widetilde {\mathcal H}(\p,U)$. 
Then one can argue as in the proof of 
\cite[Proposition 4.5.3]{DG*} (see also \cite[Proposition 5.6.4]{DG*}). \qed
 
\subsubsection{Entropy in  locally linearly compact vector spaces}
In the sequel we consider topological vector spaces over a discrete field $K$, endowed with linear topologies, i.e., having a local base of 
\nbd s of 0 formed by open subspaces. Such a space $V$ is said to be {\em linearly compact}, if every filter base in $V$ formed by 
cosets of  closed subspaces of $V$ has a non-empty intersection. A topological vector space space $V$ is said to be {\em locally linearly compact}, if 
$V$ has an open  linearly compact subspace. These notions were introduced by Lefschetz in 1942. 

In analogy to the situation observed in locally compact groups, when $V$ is a locally linearly compact vector space and $\p$ is a continuous vector space endomorphism of $V$, every open  linearly compact subspace $U$ of $V$ is  $\dim$-$\p$-inert in the sense of \S \ref{i-i-entropy} (i.e., 
 $\dim (\p(U)/(U\cap \p(U))) < \infty$). The argument to see this is similar to that of Lemma \ref{laaast:lemma}: the subspace $ \p(U)$ is linearly compact, as a continuous image of the linearly compact subspace $U$, whereas $U\cap \p(U)$ is an open subspace of the the linearly compact subspace $\p(U)$. This entails finiteness of the dimension of the (discrete) quotient $\p(U)/(U\cap \p(U))$. Therefore, $U$ is $\dim$-$\p$-inert, so one can define the intrinsic $\dim$-entropy of $\p$. This entropy was deeply investigated in \cite{CGB} (see also  \cite{CI,CGB1} for its counterpart -- the topological entropy). 
%%%%%%%%%%%%%%%%%%%%%%%%%%%%

\subsection{Epilogue} Here we only mention very briefly three more relevant topics that were not touched so far.  Before doing that let us mention that the algebraic entropy, so far defined for endomorphisms of abelian groups, was extended recently also to amenable actions of semigroups in \cite{DFG} (a semigroup $S$ is  {\em right amenable} if for every finite subset $K\subseteq S$ and every $\varepsilon>0$ there exists finite non-empty 
subset $F\subseteq S$, such that $|Fx\setminus F|\leq \varepsilon|F|$ for every $x\in K$). 
More precisely, if $S$ is a right amenable cancellative semigroup with a left action $S\overset{\alpha}{\curvearrowright}A$ of $S$ on an abelian group $A$, then the algebraic entropy $h_{alg}(S\overset{\alpha}{\curvearrowright}A)$ is defined in such a way that if $f\in\End(A)$ and the action  $\N\overset{\alpha}{\curvearrowright}A$ is defined by $\alpha(n)(x) = f^n(x)$ for $n\in \N$ and $x\in A$,  then $h_{alg}(\N\overset{\alpha}{\curvearrowright}A)$ coincides with the algebraic entropy $h_{alg}(f)$ already defined. (The definition uses a generalization of Fekete Lemma from  \cite{CCK}.) This entropy shares many of the properties of the algebriac entropy defined for single endomorphisms, for more detail see \cite{DFG} (see also \cite{V2} for algebraic entropy of actions of amenable groups on locally compact abelian groups). 

\subsubsection{Adjoint entropy and intrinsic adjoint entropy}

\def\aent{\ent^\star}

In analogy to the algebraic entropy $\mbox{ent}$, in \cite{DGS} the adjoint algebraic entropy of endomorphisms of abelian groups $A$ was introduced ``replacing" the family of all finite subgroups of $G$  by the family of all finite-index subgroups of $A$ as follows.  For a  finite-index subgroup $H$ of $A$, an endomorphism $\p:A\to A$ and $n\in\N_+$ the subgroup $C_n(\p,H)=H\cap\p^{-1}(H)\cap\ldots\cap\p^{-n+1}(H)$ still has a finite index. We call it \emph{$n$-th $\p$-cotrajectory of $H$} and we call \emph{$\p$-cotrajectory of $H$} the subgroup $C(\p,H)=\bigcap_{n\in\N}\p^{-n}(H).$ Clearly, it is the maximum $\p$-invariant subgroup of $H$.

\begin{itemize}
\item[(a)] The \emph{adjoint algebraic entropy of $\p$ with respect to $H$} is 
\begin{equation}\label{H*}
\mathcal H^\star(\p,H)={\lim_{n\to \infty}\frac{\log|A:C_n(\p,H)|}{n}}.
\end{equation}
\item[(b)] The \emph{adjoint algebraic entropy of $\p$} is 
$$
\aent(\p)=\sup\{\mathcal H^\star(\p,H):H \mbox{ is finite-index subgroup of }A)\}.
$$
\end{itemize}

The adjoint algebraic entropy, although a very natural counterpart of the algebraic one, does not have equally nice properties for the  simple reason that many groups fail to   have sufficiently many finite-index subgroups. Consequently, this entropy is not monotone with respect to taking restrictions to invariant subgroups. 

The remedy is provided again by the inert subgroups. A careful analysis of the limit in (a) shows that the sequence of the sizes of the finite quotient groups $\frac{C_n(\p,H)}{C_{n+1}(\p,H)} $ is stationary \cite{DGS}. So, in order to define the limit one only needs the finiteness of the quotients. So similarly to the case of the intrinsic entropy, it is enough to have a $\p$-inert subgroup $H$ in order to be able to define {\em intrinsic adjoint entropy} $\widetilde{\mathcal H}^\star(\p,H)$ of $\p$ with respect to $H$ by $\widetilde{\mathcal H}^\star(\p,H)={\lim_{n\to \infty}\frac{\log[H:C_n(\p,H)]}{n}}$. The {\em intrinsic adjoint entropy} of $\p$ is 
$$
\widetilde{\aent}(\p)= \sup\{\widetilde{\mathcal H}^\star(\p,H): H\in \mathcal I_\p(A)\}.
$$
This entropy has better properties than the adjoint entropy \cite{GBS1}.  More details on the adjoint entropy and the intrinsic adjoint entropy can be found in \cite{DGS,GBS1}.

\subsubsection{Entropy in normed semigroups and a unified approach to all entropies} In the expository paper \cite{DG:PC}  a unifying approach for many known entropies in many branches in mathematics, algebra, topology, measure theory, set theory, etc. is described.

The basis step of this approach is a notion of semigroup entropy $h_\mathfrak S$ in the category $\mathfrak S$ of normed semigroups and contractive homomorphisms. 
It has a number of natural properties (as the Logarithmic Law, Invariance under Conjugation, etc.).  The second step is to define, for a specific category $\mathfrak  X$, a functor $ F:\mathfrak  X \to \mathfrak S$ and then defined the entropy $h_F$ as the composition $h_F=h_\mathfrak S\cdot F$. More precisely, for every $\mathfrak X$ morphism $\p$ one takes $h_F:= h_\mathfrak S(F\p)$. This provides, under certain 
natural restrains on $F$, the same (or the counterparts of the) properties proved for $h_\mathfrak S$. This general scheme permits to obtain many of the known entropies as $h_F$, for appropriately chosen categories $\mathfrak  X$ and functors $ F:\mathfrak  X \to \mathfrak S$. The fundamental properties of the algebraic entropy for group endomorphisms easily follow using this setting. 

\subsubsection{Topological entropy and its relation to algebraic entropy}  In \S 6.5 we defined the algebraic entropy $h_{alg}$ for continuous endomorphisms of a locally compact group $G$, endowed with a right Haar measure $\mu$.  In what follows we discuss the {\em topological entropy} $h_{top}$ for such endomorphisms.  

 For $n\in\N_+$ and $U\in \mathcal U_G$, the {\em $n$-th $\psi$-cotrajectory} $C_n(\psi,U)=U\cap \psi^{-1}(U)\cap \dots\cap \psi^{-n+1}(U)$  is compact, so measurable. Let 
\begin{equation}\label{h_mu(U)} 
\mathcal 
H_{top}(\psi,U)=\limsup_{n\to\infty}-\frac{\log\mu(C_n(\psi,U))}{n}.
\end{equation}
The topological entropy is defined by 
\begin{equation}\label{h_mu} 
h_{top}(\psi)=\sup\{\mathcal H_{top}(\psi,U):U\in \mathcal U_G\}
\end{equation} 
 (see \cite{DG*,DSV}).
Moreover, if $G$ totally disconnected locally compact group and $U\in\mathcal B(G)$, then  $C_n(\psi,U)$ is an open subgroup of the compact group $U$, hence $|G: C_n(\psi,U)|< \infty$ and one can obtain the following measure-free version of the formula (\ref{h_mu(U)})
$$
\mathcal 
H_{top}(\psi,U)=\lim_{n\to\infty}\frac{\log |U: C_n(\psi,U)|}{n}.
$$
Finally, if $G$ is also compact (i.e., profinite), then one can take $U=G$ obtaining directly this convenient counterpart of (\ref{h_mu(U)})
$$
H_{top}(\psi,U)=\lim_{n\to\infty}\frac{\log |G: C_n(\psi,U)|}{n}.
$$
The utmost interest in  topological entropy stems from topological dynamics, this is why the relevance of  the algebraic entropy
is motivated by its close connection to the topological one which we discuss now. 
Following \cite{DG_BT}, we call {\em Bridge Theorem} statements of the form $h_{top}(\p)= h_{alg}(\widehat \p)$ 
for a continuous endomorphism $\p:G\to G$ of a locally compact abelian group and its Pontryagin dual $\widehat \p: \widehat G \to \widehat G$. 
The Bridge Theorem was proved in \cite{DG_BT} for $G$ compact, and previously by Weiss \cite{We} when $G$ is profinite and by 
Peters \cite{Pet}, when $G$ is compact metrizable and $\p$ is a topological automorphism. In \cite{DG_BT1} the Bridge Theorem is proved for totally disconnected locally compact abelian groups. A proof of the Bridge theorem for topological automorphisms  of locally compact abelian groups can be found in \cite{V2}.

\subsubsection{Entropy and scale function}
%%%%%%%%%%%%%%%%%%%%%%%%%%%
The scale function was introduced by Willis \cite{Willis1,Willis2} as follows.

Let $G$ be a locally compact totally disconnected group.  According to Lemma \ref{laaast:lemma},  every $U\in \mathcal B_G$ is $\p$-inert for every $\p\in C\End(G)$.  This means that $s_G(\p,U):= |U^\p:(U^\p\cap U)|$ is a positive natural number. The \emph{scale} of $\phi$ is 
\begin{equation}\label{(S)}
s_G(\phi) = \min\{s_G(\phi,U): U \in \mathcal B_G\}. 
\end{equation}
We use the notations $s(\phi,U)$ and $s(\phi)$ when the group $G$ is clear from the context.

Since the scale function is defined as a minimum, there exists $U \in \mathcal B_G$ for which this minimum realizes, that is $s(\phi) =s(\phi,U)$, and such $U$ is called \emph{minimizing} for $\phi$ in \cite{Willis2}. Let $\mathcal M(G,\phi)$ be the subfamily of $\mathcal B_G$ consisting of all minimizing subgroups  for $\phi$.
  In  \cite{Willis1} the scale function was defined 
 only for internal automorphisms of $G$ and its definition was motivated by the fact that in the non-commutative case, 
 a locally compact totally disconnected group $G$ may fail to have sufficiently many compact open normal subgroups, 
 in particular, $\mathcal B_G$ may contain no normal subgroups at all, when $G$ is neither compact nor discrete. 

\smallskip
It is easy to see that $\mathcal M(G,\phi)$ contains all $\phi$-invariant;
in particular,  $s(\phi)=1$ precisely when there exists a $\phi$-invariant $U\in\mathcal B_G$. In particular, $s(\phi)=1$
for all $\p$ when $G$ is compact (as $s(\phi,G)=1$).

For a topological automorphism $\phi:G\to G$ of  a totally disconnected locally compact group $G$, the subgroup
\begin{equation}\label{DefNUB}
\nub(\phi)=\bigcap\{U\in\mathcal{B}_G: U\in\mathcal M(G,\phi)\}.
\end{equation}
 is $\phi$-stable and compact.

The following theorem, proved in \cite{BDGB}, answered a question posed by T. Weigel in 2011, 
about the relation between the scale function and the entropy.
 
\begin{thm}\label{nubbanale}
Let $G$ be a totally disconnected locally compact group and $\phi:G\to G$ a topological automorphism. Then $\log s(\phi)\leq h_{top}(\phi)$, the equality 
$\log s(\phi)=h_{top}(\phi)$ holds if and only if $\nub(\phi)=\{e_G\}$.
\end{thm}

 For the relation of the topological entropy with the scale see also  the recent paper \cite{GV3}.

%------------------------------------------------------------------------------------%

%\bigskip

%\vskip -0.4 true cm

\begin{center}{\textbf{Acknowledgments}}
\end{center}
The authors wish to thank Alma D'Aniello for her useful comments during the preparation of the paper.  
\vskip -0.5 true cm

%------------------------------------------------------------------------------------%

%%%%%%%%%BIBLIOGAFIA

%%%%%%%%%%%%

%%%%%%%%%%%%%%%%%%%%%%%%%%%%%%%%%%%%%

%\newpage
%\vskip -0.5mm

\def\pn{\par\noindent}
\def\cen{\centerline}

{\footnotesize \pn{\bf Ulderico Dardano}\; \\ {Dipartimento di Matematica e
Applicazioni ``R.Caccioppoli'}, {Universit\`a di Napoli ``Federico
II'', Via Cintia - Monte S. Angelo, I-80126,} {Napoli, Italy}\\
{\tt Email: dardano@unina.it}
\medskip

{\footnotesize \pn{\bf Dikran Dikranjan}\; \\ {Dipartimento di Matematica ed Informatica}, {Universit\`a di Udine, Via delle Scienze, 208, Udine, Italy, I-33100}\\
{\tt Email: dikranjan@dimi.uniud.it}
\medskip 

{\footnotesize \pn{\bf Silvana Rinauro}\; \\ {Dipartimento di Matematica, Informatica ed Economia,} {Universit\`a della
Basilicata, Viale dell'Ateneo Lucano 10,
I-85100,} { Potenza, Italy}\\
{\tt Email: silvana.rinauro@unibas.it}\\ 
%%%%%%%%%%%%%%%%%%%%%%%%%%%%%%%%%%%%%%%%

\end{document}